\documentclass[9pt,final]{IEEEtran}
\title{Optimal Control with Noisy Time}
\author{Andrew Lamperski\thanks{Department of Engineering, University
    of Cambridge, Cambridge, UK ({\tt a.lamperski@eng.cam.ac.uk}).}
  \and Noah J. Cowan\thanks{Department of Mechanical Engineering, The
    Johns Hopkins University, Baltimore, MD, USA ({\tt
      ncowan@jhu.edu}).}}

\usepackage[dvipsnames]{color}

\usepackage{cite}
\usepackage{amsmath,amsfonts,amssymb}
\usepackage{mathtools,enumerate}
\usepackage{graphicx}
\graphicspath{{figures/}}

\newtheorem{theorem}{Theorem}
\newtheorem{lemma}{Lemma}
\newtheorem{problem}{Problem}
\newtheorem{proposition}{Proposition}

\newtheorem{remark}{Remark}
\newtheorem{example}{Example}

\def\tp{\mathsf{T}}
\newcommand{\ind}{\mathbf{1}}
\newcommand{\Tr}{\mathrm{Tr}}
\renewcommand{\P}{\mathbb{P}}
\newcommand{\F}{\mathcal{F}}
\newcommand{\A}{\mathcal{A}}
\newcommand{\R}{\mathbb{R}}
\newcommand{\C}{\mathbb{C}}
\newcommand{\N}{\mathbb{N}}
\newcommand{\B}{\mathcal{B}}
\newcommand{\D}{\mathcal{D}}
\newcommand{\X}{\mathcal{X}}
\newcommand{\U}{\mathcal{U}}
\newcommand{\Cont}{\mathcal{C}}
\newcommand{\E}{\mathbb{E}}
\newcommand{\dom}{\mathrm{dom}}
\newcommand{\spec}{\mathrm{spec}}
\renewcommand{\vec}{\mathrm{vec}}
\newcommand{\gap}{\vspace{.3cm}}
\renewcommand{\Re}{\mathrm{Re}\:}
\DeclareMathOperator*{\argmin}{arg\,min}

\begin{document}
\maketitle

\begin{abstract}
  This paper examines stochastic optimal control problems in which the
  state is perfectly known, but the controller's measure of time is a
  stochastic process derived from a strictly increasing L\'evy
  process.  We provide dynamic programming results for continuous-time
  finite-horizon control and specialize these results to solve a
  noisy-time variant of the linear quadratic regulator problem and a
  portfolio optimization problem with random trade activity rates. For
  the linear quadratic case, the optimal controller is linear and can
  be computed from a generalization of the classical Riccati
  differential equation.
\end{abstract}

\section{Introduction}

Effective feedback control often requires accurate timekeeping. For
example, finite-horizon optimal control problems generally result in
policies that are time-varying functions of the state.
However, chronometry is imperfect and thus feedback laws are
inevitably applied at incorrect times. Little appears to be known
about the consequences of imperfect timing on control
\cite{lavalletime2007,carverstateestimation2013,lamperskitimechanged2013}.
This paper addresses optimal control with temporal uncertainty.

A stochastic process can be \emph{time-changed} by replacing its time
index by a monotonically increasing stochastic process
\cite{veraarttime2010}.
Time-changed stochastic processes arise in finance, since changing the
time index to a measure of economically relevant events, such as
trades, can improve modeling
\cite{clarksubordinated1973,aneorder2000,carrtimechanged2004}.  This
new time index is, however, stochastic with respect to ``calendar''
time.

We suspect that similar notions of stochastic time changing may
facilitate the study of time estimation and movement control in the
nervous system. Biological timing is subject to noise and
environmental perturbation \cite{eaglemanhuman2008}. Furthermore,
humans rationally exploit the statistics of their temporal noise
during simple timed movements, such as button pushing
\cite{jazayeritemporal2010} and pointing \cite{hudsonoptimal2008}.  To
analyze more complex movements, a theory of feedback control that
compensates for temporal noise seems desirable.

Within control, the most closely related work to the present paper
deals with analysis and synthesis of systems with uncertain sampling
times.  The study of uncertain sampling times has a long history in
control \cite{kushnerstability1969}, and is often motivated by
problems of clock jitter \cite{wittenmarktiming1995,skafanalysis2009}
or network delays \cite{hespanhasurvey2007}. In these works, control
inputs are sampled at known times and held over unknown intervals. To
derive the dynamic programming principle in this paper, system
behavior is analyzed for control inputs held over random intervals,
bearing some similarity to optimal control with random sampling
\cite{adesstochastic2000}. Fundamentally, however, studies of sampling
uncertainty assumes that an accurate clock can measure the sample
times; the present work relaxes this assumption.

Other aspects of imperfect timing have been addressed in control
research to a more limited extent.  For example, the importance of
synchronizing clocks in distributed systems seems clear
\cite{frerisfundamental2007,simeonedistributed2008}, but more work is
needed to understand the the implications of asynchronous clock
behavior on common control issues, such as stability
\cite{lorandstability2003} and optimal performance
\cite{singhlqr2011}. 


This paper focuses on continuous-time stochastic optimal control 
with
perfect state information, but a stochastically time-changed control
process. 
Dynamic programming principles for general nonlinear stochastic
control problems are derived, based on extensions of the classical
Hamilton-Jacobi-Bellman equation.  The results apply to a wide class
of stochastic time changes given by strictly increasing L\'evy
processes.  The dynamic programming principles are then specialized to
give explicit solutions to time-changed versions of the finite-horizon
linear quadratic regulator and a portfolio optimization problem.

Section~\ref{sec:prelim} defines the notation used in the paper,
states the necessary facts about L\'evy,
and defines the class of noisy clock models used. The main results on
time-changed diffusions and optimal control are given in
Section~\ref{sec:results}. The results are proved in
Sections~\ref{sec:pf} with supplementary arguments given in the
appendices. Sections~\ref{sec:discussion} and~\ref{sec:conclusion} discuss future work and conclusions, respectively. 

\section{Preliminaries}\label{sec:prelim}

After establishing notation and reviewing L\'evy processes, this
section culminates in the construction of L\'evy-process-based clock
models upon which the remainder of the theory of this paper is built.

\subsection{Notation}\label{sec:notation}
The norm symbol, $\|\cdot\|$, is used to denote the Euclidean norm for
vectors and the Frobenius norm for matrices. 

For a set $S$, its closure is denoted by $\overline{S}$. 

The spectrum of matrix $A$ is denoted by $\spec(A)$. 

The Kronecker product is denoted by $\otimes$, while the Kronecker sum
is denoted by $\oplus$:
\begin{equation*}
A\oplus B = A\otimes I + I\otimes B.
\end{equation*}
The vectorization operation of stacking the columns of a matrix is
denoted by $\vec$.

A function $h:\R\times\R^n\to \R$ is in $\Cont^{1,2}$ if $h(s,x)$ is continuously differentiable in $s$, twice continuously differentiable in $x$. The function $h$ is said to satisfy a \emph{polynomial growth condition}, if in addition, there are constants $K$ and $q$ such that 
\begin{equation*}
|h(s,x)|,\:\left|
  \frac{
    \partial h(s,x)}{\partial s} 
\right|,
\:
\left|
  \frac{\partial h(s,x)}{\partial x_i}
\right|,\:
\left|
  \frac{\partial^2 h(s,x)}{\partial x_i\partial x_j}
\right|
\le K\left(1+\|x\|^q\right),
\end{equation*}
for $i,j = 1,\ldots n$, and all $x\in \R^n$. In this case, $h\in\Cont_p^{1,2}$ is written.  

Stochastic processes will be denoted as $\zeta_t$, $X_s$, etc., with time
indices as subscripts. Occasionally,
processes with nested subscripts will be written with parentheses,
e.g. $\zeta_{\tau_s} = \zeta(\tau_s)$. Similarly, the elements of a
stochastic vector will be denoted as $X_1(s)$. 

Functions that are right-continuous with left-limits will be called \emph{c\`adl\`ag}, while functions that are left-continuous with right-limits will be called \emph{c\`agl\`ad}.

\subsection{Background on L\'evy Processes}

Basic notions from L\'evy processes required to define the general
class of clock models are now reviewed. The definitions and results
can be found in \cite{applebaumlevy2004}.

A real-valued stochastic process $Z_s$ is called a \emph{L\'evy process} if 
\begin{itemize}
\item $Z_0=0$ almost surely (a.s.).
\item $Z_s$ has independent, stationary increments: If $0\le r\le s$, then $Z_r$ and $Z_s-Z_r$ are independent and
  $Z_s-Z_r$ has the same distribution as $Z_{s-r}$.
\item $Z_s$ is stochastically continuous: For all $a>0$ and all $r\ge
  0$, $\lim_{r\to s}\P(|Z_s-Z_r|>a) = 0$. 
\end{itemize}

It will be assumed that L\'evy processes in this paper are
right-continuous with left-sided limits, i.e. they are
c\`adl\`ag.  No generality is lost since, for every L\'evy
process, $Z_t$, there is a c\`adl\`ag L\'evy process, $\tilde{Z}_t$,
such that $Z_t = \tilde Z_t$ for almost all $t$.

Some of the more technical arguments rely on the notion of Poisson
random measures, which will now be defined. Let $\B$ be the Borel subsets of $\R$ and let $(\Omega,\Sigma,\P)$ be a probability space. A \emph{Poisson random
  measure} is a function $N:[0,\infty) \times \B\times \Omega \to \N\cup\{\infty\}$, such that
\begin{itemize}
\item For all $s\ge 0$ and $\omega\in\Omega$, $N(s,\cdot,\omega)$ is a measure.
\item For all disjoint Borel subsets $A,B\in \B$ such that $0\notin
  \overline{A}$ and $0\notin\overline{B}$, $N(\cdot,A,\cdot)$ and $N(\cdot,B,\cdot)$ are
  independent Poisson processes.
\end{itemize}
Typically, the $\omega$ argument will be dropped, and it will be
implicitly understood that $N(s,A)$ denotes a measure-valued
stochastic process. 

The following relationship between L\'evy processes and Poisson random
measures will be used in several arguments. For a L\'evy process,
$Z_s$, with jumps denoted by 
$\Delta Z_s$, there is a Poisson random measure that counts the number
of jumps into each Borel set $A$ with $0\notin \overline{A}$:
\begin{equation*}
N(s,A) = \left|\{\Delta Z_r \in A : 0\le r \le s \}\right|.
\end{equation*} 


\textbf{Subordinators.} A monotonically increasing L\'evy process, $\tau_s$, is called a
\emph{subordinator}. 
The following properties of subordinators will be
used throughout the paper.
\begin{itemize}
\item \textbf{Laplace Exponent:} There is function, $\psi$, called the {\it
    Laplace exponent}, defined
  by
  \begin{equation}\label{laplaceForm}
    \psi(z) = bz + \int_0^{\infty}\left(1-e^{-zt}\right)\lambda(dt),
  \end{equation}
  such that 
  \begin{equation}\label{eq:laplace}
    \E\left[
      e^{-z\tau_s}
    \right] = e^{-s\psi(z)} \qquad \textrm{for all} \qquad z\ge 0.
  \end{equation}
  Here $b\ge 0$ and the measure satisfies
  $\int_0^{\infty}\min\{t,1\}\lambda(dt) < \infty$. The measure $\lambda$ is called a \emph{L\'evy measure}. The pair
  $(b,\lambda)$ is called the \emph{characteristics} of $\tau_s$.
\item \textbf{L\'evy-It\^o Decomposition:} There is a Poisson random
  measure $N$ such that 
  \begin{equation*}
  \tau_s = bs + \int_0^{\infty}t N(s,dt).
  \end{equation*}  
  Furthermore, if $A\subset (0,\infty)$ is a Borel set such that
  $0\notin \overline{A}$, then $\E[N(1,A)] = \lambda(A)$.
\end{itemize}





\gap

The function, $\psi$, is called the Laplace exponent because
\eqref{eq:laplace} is the Laplace transform of the distribution of
$\tau_s$. 

For control problems, simpler formulas will often result from
replacing $\psi$ with the function $\beta(z) = -\psi(-z)$. Note then,
that $\beta$ has the form 
\begin{equation}\label{betaDef}
  \beta(z) = bz + \int_0^{\infty}\left(
    e^{zt}-1
  \right)\lambda(dt).
\end{equation}

Define $r_{\max}$ by 
\begin{equation*}
r_{\max} = \sup\left\{r:
  \int_1^{\infty} e^{rt} \lambda(dt)
\right\}
\end{equation*}
and define the domain of $\beta$ as
\begin{equation*}
\dom(\beta) = \{z\in\C: \Re z < r_{\max}
\}.
\end{equation*}
Note that $\int_1^{\infty} \lambda(dt) < \infty$ implies that $r_{\max}\in [0,\infty]$. 

The function $\beta$ is used to construct optimal solutions for the
linear quadratic problem, as well as the portfolio problem below. The
main properties are given in the following lemma, which is proved in Appendix~\ref{pf:beta}.

\gap

\begin{lemma}\label{lem:beta}{\it
  For all $z\in \dom(\beta)$, the function $\beta$ is analytic at $z$,
  and 
  \begin{equation}\label{betaExp}
    \E\left[
      e^{z\tau_s}
    \right] = 
    e^{s\beta(z)}.
  \end{equation}
  Furthermore, if $A$ is a square matrix with $\spec(A)\subset\dom\beta$, then
  \begin{equation}\label{matrixBeta}
    \beta(A) = b A + \int_0^{\infty}\left(e^{At}-I\right) \lambda(dt)
  \end{equation}
  is well defined and 
  \begin{equation}
    \E\left[e^{A\tau_s}\right] = e^{s\beta(A)}. \label{betaM}
  \end{equation}
}
\end{lemma}

\gap

Since $\beta$ is analytic, several
methods exist for numerically computing the matrices $\beta(A)$ \cite{highamfunctions2008}. In some special cases, as discussed below, $\beta(A)$ may be computed using well-known matrix computation methods. 

\gap

\begin{example}{\rm \label{ex:Poisson}
    The simplest non-trivial subordinator is the Poisson process $N_t$,
    which is characterized by 
    \begin{equation*}
    \P(N_t = k) = e^{-\gamma t} \frac{(\gamma t)^k}{k!},
    \end{equation*}
    where $\gamma >0$ is called the rate constant. Its Laplace exponent is
    given by $\psi(z) = \gamma - \gamma e^{-z}$, which is found by
    computing the expected value directly. The characteristics are
    $(0,\gamma\delta(t-1))$. In this case, $\dom(\beta) = \C$, and $\beta(A) = \gamma e^A - \gamma
    I$, which can be computed from the matrix exponential.
  }
\end{example}

\gap

\begin{example}{\rm
    The gamma subordinator, which is often used to model ``business time''
    in finance \cite{cvitanicoptimal2008,madanvariance1998}, has
    increments distributed as gamma random 
    variables. It has Laplace exponent $\psi(z) = \delta\log(1+z/\gamma)$
    with characteristics $b=0$ and $\lambda(dt) = \delta e^{-\gamma
      t}t^{-1} dt$. Thus $\beta(z) = -\delta\log(1-z/\gamma)$,
    $\dom(\beta) = \{z\in\C: \Re z < \gamma\}$, and matrix function 
    $
    \beta(A) = -\delta \log\left(I-\gamma^{-1} A\right)
    $
    may be computed from the matrix logarithm.
  }
\end{example}

\gap

\textbf{Why L\'evy Processes?} In the next subsection, the clock model
in this paper will be constructed from a subordinator $\tau_s$. The
motivation for using
L\'evy processes will be explained. Consider a continuous-time noisy clock, $c_s$ which
is sampled with period
$\delta$. A natural model might take the form
\begin{equation}\label{dtClock}
  c_{\delta (k+1)} = c_{\delta k}+ \delta + n(k,\delta),
\end{equation}
where $n(j,\delta)$ are random variables. In this case, the clock
increments consist of a deterministic step of magnitude $\delta$ plus
a random term. 

If $c_s$ is a L\'evy process, then by definition, all of the
increments $c_{\delta (k+1)}-c_{\delta k}$ are independent and
identically distributed. Thus, the decomposition in \eqref{dtClock} holds with
$n(k,\delta) = c_{\delta (k+1)} - c_{\delta k} - \delta$. 
If $c_s$
were not a L\'evy process, then \eqref{dtClock} may hold for some
particular $\delta$, but there might be another period,
$\delta' < \delta$, for which the decomposition fails. The L\'evy
process assumption will guarantee that the clocks are well-behaved
when taking continuous time limits (i.e. $\delta \downarrow 0$).

\subsection{Clock Models}\label{sec:subjectiveTime}

 
Throughout the paper, $t$ will denote the time index of the
plant dynamics, while $s$ will denote the value of clock available to
the controller. Often, $t$ and $s$ will be called {\it plant time} and {\it controller time}, respectively.  
The interpretation of $s$ and $t$ varies depending on context.  In
biological motor control, $t$ would denote real time, since the limbs
obey Newtonian mechanics with respect to real-time, while $s$ would
denote the internal representation of time.  For the portfolio problem
studied in Subsection~\ref{sec:dp}, an opposite interpretation
holds. Here, the controller (an investor) can accurately measure
calendar time, but price dynamics are simpler with respect a different
index, ``business time'', which represents the progression of economic
events
\cite{aneorder2000,carrtimechanged2004,clarksubordinated1973}. Thus,
$s$ would denote calendar time, while $t$ would denote business time,
which might not be observable.

The relationship between $s$ and $t$ will be described
stochastically. 
Let $\tau_s$ be a strictly increasing subordinator. In other words, if
$s<s'$ then $\tau_{s}<\tau_{s'}$ a.s. 
(Note that any subordinator can be
made to be strictly increasing by adding a drift term $bs$ with $b>0$.) The
process $\tau_s$ will have the 
interpretation of being the amount of plant time that has passed
when the controller has measured $s$
units of time. The process $\zeta_t$ will be an
inverse process that describes how much time the controller measures
over  $t$ units of plant time. Formally, $\zeta_t$ is defined by
\begin{equation}\label{zetaDef}
  \zeta_t = \inf \{\sigma: \tau_{\sigma} \ge t\}.
\end{equation}
Note that $\zeta(\tau_s) = s$ a.s. Indeed, $\zeta(\tau_s) =
\inf\{\sigma:\tau_{\sigma} = \tau_s\}$, by definition. Since $\tau_s$ is right
continuous and strictly
increasing, a.s., it follows that $\zeta_{\tau_s} = s$, a.s.

\gap

\begin{example}{\rm
    The case of no temporal uncertainty corresponds to $\tau_s = s$ and
    $\zeta_t = t$. The Laplace exponent of $\tau_s$ is computed directly as $\psi(z) =
    z$ and the characteristics are $(1,0)$. Here $\dom(\beta) = \C$. 
  }
\end{example}

\gap

\begin{example}\label{ex:IG}{\rm A more interesting temporal noise
    model, also used as a ``business time'' model
    \cite{barndorffnielsenprocesses1998}, is the inverse Gaussian
    subordinator. Fix $\gamma>0$ and $\delta>0$. Let $C_t = \gamma
    t+W_t$, where $W_t$ is a standard unit Brownian motion. The
    inverse Gaussian subordinator is given by \begin{equation*} \tau_s
      = \inf\{t:C_t=\delta s\}, \end{equation*} with Laplace exponent
    $\psi(z) = \delta(\sqrt{\gamma^2+2z}-\gamma)$. Here $b=0$ and
    $\lambda$ is given by \begin{equation*} \lambda(dt) =
      \frac{\delta}{\sqrt{2} \Gamma(1/2)} e^{-\frac{1}{2}\gamma^2
        t}t^{-\frac{3}{2}} dt , \end{equation*} where $\Gamma$ is the
    gamma
    function.  
    Here, $\dom(\beta)$ corresponds to $\Re z < \gamma^2 /2$ and
    $\beta(A) = \delta\left(\gamma I - \sqrt{\gamma^2 I - 2
        A}\right)$, which can be computed from the matrix square root.
    It can be shown that the inverse process is given
    by \begin{equation*} \zeta_t =
      \sup\left\{\delta^{-1}C_{\sigma}:0\le\sigma\le t
      \right\}.  \end{equation*} See Figure
    \ref{fig:IGnoise}.  

    \begin{figure}
      \centering
      \includegraphics{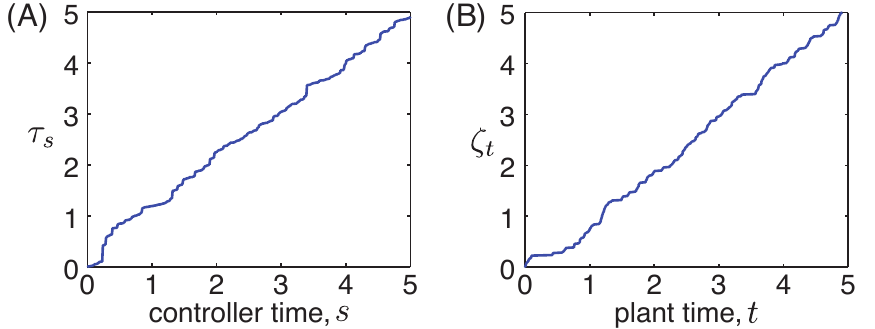}
      \caption{\label{fig:IGnoise} (A). The inverse Gaussian
        subordinator, $\tau_s$, with $\gamma = \delta =2$. The process
        was simulated by generating independent inverse Gaussians
        using the method from \cite{michaelgenerating1976}.  (B) The
        inverse process, $\zeta_t$. Note that the graph of $\zeta_t$
        can be found from the graph of $\tau_s$ by simply switching
        the axes.  }
    \end{figure}
  }
\end{example}

\gap

In the preceding example, the process $\tau_s$ has jumps,
but the inverse, $\zeta_t$, is continuous. The next proposition
generalizes this observation for any strictly increasing
subordinator, $\tau_s$. 

\gap

\begin{proposition}\label{prop:cont}
  {\it
    The process $\zeta_t$ is continuous almost surely.
  }
\end{proposition}

\gap

\begin{IEEEproof}
  Fix $\epsilon>0$ and  $t\ge 0$. Set $s=\zeta_t$. Strict
  monotonicity of $\tau_s$ implies that
  $[\tau_{\max\{s-\epsilon,0\}},\tau_{s+\epsilon}]$ is a nonempty
  interval, a.s. The inverse property of $\zeta_t$ implies (almost surely) that $t\in
  [\tau_{\max\{s-\epsilon,0\}},\tau_{s+\epsilon}]$ and $\zeta_{t'}\in
  [\max\{s-\epsilon,0\},s+\epsilon]$ for all $t'\in
  [\tau_{\max\{s-\epsilon,0\}},\tau_{s+\epsilon}]$. 
\end{IEEEproof}

\section{Main Results}\label{sec:results}
This section presents the main results of the paper. First, given an
It\^o process, $Y_t$, a representation of the time-changed process $X_s
= Y(\tau_s)$ as a semimartingale with respect to controller time, $s$, is
derived. This representation is then used to derive a general dynamic
programming principle for control problems with noisy clocks. As an
example, the dynamic programming principle is used to solve a simple
portfolio optimization problem under random trade activity
rates. Finally, the dynamic programming method is used to solve a
noisy-time variant of the linear quadratic regulator problem. All
proofs are given in Section~\ref{sec:pf}.    

\subsection{Time-Changed Stochastic Processes}
\label{sec:timeChange}

This section gives a basic representation theorem for time-changed
stochastic processes that will be vital for dynamic programming
proofs. The theorem is proved in Subsection~\ref{pf:tcDiffusion}. 

Let $W_t$ be a Brownian motion with $\E[W_tW_t^\tp] = tI$. Let $Y$ be a stochastic process 
defined by
\begin{equation}\label{diffusion}
  dY_t = F_tdt + G_tdW_t,
\end{equation}
where $F_t$ and $G_t$ are $\F_t^W$ predictable processes, where
$(\F_t^W)_{t\ge 0}$ is the $\sigma$-algebra generated by
$W_t$. Furthermore, assume that $F_t$ and $G_t$ are left-continuous
with right-sided limits.  


Let $\F^{\tau,W}=(\F_s^{\tau,W})_{s\ge 0}$ be the smallest filtration
such that for all $r\in [0,s]$ and all $t\in [0,\tau_s]$ both $\tau_r$
and $W_t$ are measurable.  

\begin{theorem}
  \label{timeChangedDiffusion}
  {\it
    Let $\tau_s$ be a subordinator characterized by $(b,\lambda)$. If the terms of \eqref{diffusion} satisfy
    \begin{itemize}
    \item $\int_0^{\tau_S}\|F_t\|dt <\infty$ almost surely and
    \item $\E\left[\int_0^{\tau_S}\|G_t\|^2dt\right] < \infty $,
    \end{itemize}
    then
    the time-changed process $X_s = Y(\tau_s)$ is an $\F^{\tau,W}$ semimartingale given by 
    \begin{multline}\label{timeChange}
      X_s = X_0 + b\int_0^sF(\tau_{r^-})dr +
      \sqrt{b}\int_0^sG(\tau_{r^-})d\tilde{W}_{r} + \\
      \sum_{0\le r \le s}\left(
        \int_{\tau_{r^-}}^{\tau_r}F_tdt + 
        \int_{\tau_{r^-}}^{\tau_r}G_tdW_t
      \right).
    \end{multline}
    Here $\tilde{W}_s$ is an $\F^{\tau,W}$-measurable Brownian motion
    defined by 
    $$\sqrt{b}\tilde{W}_s = W(\tau_s) - \sum_{0\le r\le s}\left(
      W(\tau_r)-W(\tau_{r^-})
    \right),
    $$
    satisfying $b\E[\tilde{W}_s\tilde{W}_s^\tp] = bs I$.
    Furthermore, 
    \begin{enumerate}
    \item \label{finiteVariation}
      $b\int_0^{s}F_{\tau_{r^-}}dr+\sum_{0\le r \le s}\int_{\tau_{r^-}}^{\tau_r}F_tdt$  has finite variation, and 
    \item \label{martingale}
      $\sqrt{b}\int_0^{s}G_{\tau_{r^-}}d\tilde{W}_r + \sum_{0\le r \le s}\int_{\tau_{r^-}}^{\tau_r} G_tdW_t$ is an $\F^{\tau,W}$ martingale.
    \end{enumerate}
  }
\end{theorem}

\subsection{Dynamic Programming}
\label{sec:dp}
This subsection introduces the general control problem studied in
this paper. First, the basic notions of controlled time-changed
diffusions and admissible systems are defined. Then, the
finite-horizon control problem is stated, and the associated dynamic
programming verification theorem is stated. 

\textbf{Controlled Time-Changed Diffusions.}
Consider a controlled diffusion
\begin{equation}\label{NLCS}
  dY_t = F(\zeta_t,Y_{t^-},U(\zeta_t)) dt + G(\zeta_t,Y_{t^-},U(\zeta_t))dW_t, 
\end{equation}
with state $Y$ and input $U$. Recall that $\zeta_t$ is defined in \eqref{zetaDef} as the
inverse process of a subordinator, $\tau_s$. Let $X_s$ denote the
time-changed process, $X_s = Y(\tau_s)$. The processes, $X_s$ is thus
a time-changed controlled diffusion. 

\textbf{Admissible Systems}
For $s\ge 0$, let $\F_s^{\zeta,X}$ be the $\sigma$-algebra generated by $(s,X_s)$, and let $\F^{\zeta,X}$ be the associated filtration. 

Let $\X\subset \R^n$ and $\U\subset\R^p$ be a set of states and a set
of inputs, respectively. A state and input trajectory ($X_s,U_s$) is
called an \emph{admissible system} if  
\begin{itemize}
\item $X_s\in\X$ for all $s\ge 0$ 
\item $U_s$ is a c\`agl\`ad, $\F^{\zeta,X}$-adapted process such
  that $U_s\in\U$ for all $s\ge 0$. 
\end{itemize}

Note that the requirement that $U_s$ is c\`agl\`ad and
$\F^{\zeta,X}$-adapted implies that $U(\zeta_t)$ may depend on the
``noisy clock'' process, $\zeta_t$, as well as $X_r$, with
$r<\zeta_t$. If $\zeta_t\ne t$, then $U(\zeta_t)$ cannot directly
measure $t$.  

\gap

\begin{problem}\label{prob:FH}
  The time-changed optimal control problem over time horizon $[0,S]$ is to find a policy $U_s$ that solves   
  \begin{equation*}
  \min_U \E\left[
    \int_0^S c(s,X_s,U_s) ds + \Psi(X_S)
  \right],
  \end{equation*}
  where the minimum is taken over all admissible systems $(X_s,U_s)$. 
\end{problem}

\gap

Given a policy, $U$, and $(s,x)\in [0,S]\times \R^n$, the cost-to-go
function $J(s,x;U)$, is defined by
\begin{equation*}
  J(s,x;U) = 
  \E\left[
    \int_s^S c(s,X_r,U_r)dr
    +\Psi(X_S)
    \:\vline \:
    X_s = x
  \right].
\end{equation*}
Note, then, that the optimal control problem can be equivalently cast
as minimizing $J(0,x;U)$ over all admissible systems. 

\textbf{Backward Evolution Operator.} As in standard continuous-time
optimal control, the \emph{backward evolution operator},
\begin{multline}\label{generatorLim}
  \A^u h(s,x) =\\
 \lim_{\sigma\downarrow 0} \frac{1}{\sigma} \left(
    \E\left[h(s+\sigma,X_{s+\sigma})|X_s = x,\: U_r = u\right]
    -h(x)
  \right),
\end{multline}
is used to formulate the dynamic programming equations. 

To calculate an explicit form for $\A^u$, an auxiliary stochastic
process is introduced.
For  $(s,x,u)\in[0,S)\times \X\times \U$, define $Y_{st}^{xu}$ by
\begin{equation}\label{constUDrift}
  Y_{st}^{xu} = x + \int_0^t F(s,Y_{sr}^{xu},u) dr + \int_0^tG(s,Y_{sr}^{xu},u)d\hat{W}_r,
\end{equation}
where $\hat{W}_r$ is a unit Brownian motion independent of $W_t$ and $\tau_s$.

Now the domain of $\A^u$ is defined. 
Let $\D$ be the set of $h\in\Cont_p^{1,2}$ such that there exist $K$
and $q$ satisfying
\begin{equation}\label{bddCond}
  \int_0^{\infty}\left|
    \E_{\hat{W}}\left[
      h(s,Y_{st}^{xu})
    \right]-h(s,x)
  \right|\lambda(dt) < K(1+\|x\|^q+\|u\|^q)
\end{equation}
for all $(s,x,u)\in[0,S) \times\X\times\U$.

It will be shown in Subsection~\ref{pf:FH} that for $h\in\D$, the
backward evolution operator for $X_s$ is given by
\begin{multline}\label{generator}
  \A^u h(s,x) = 
  \frac{\partial h(s,x)}{\partial s} 
  +b\frac{\partial h(s,x)}{\partial x} F(s,x,u)  \\
  + \frac{1}{2}b\Tr\left(
    G(s,x,u)^{\tp} \frac{\partial^2 h(s,x)}{\partial x^2} G(s,x,u) 
  \right)
  \\
  +
  \int_0^{\infty}(\E_{\hat{W}}[h(s,Y_{st}^{xu})] - h(s,x))\lambda(dt).
\end{multline}

\gap
\begin{remark}{\rm
    When the dynamics are time-homogeneous, i.e. $F(s,y,u) =
    F(y,u)$ and $G(s,y,u) = G(y,u)$, and the policy is Markov, $U_s =
    U(X_{s^-})$, the expression for $\A^u$ in \eqref{generator} is a special case of
    Phillips' Theorem \cite{applebaumlevy2004,satolevy1999}. In this case,
    the formula can be derived using techniques from semigroup theory
    \cite{satolevy1999}. The derivation in this paper is instead based on
    It\^o calculus. 
  }
\end{remark}

\gap

\textbf{Finite Horizon Verification.} The following result is a dynamic programming verification theorem for
Problem~\ref{prob:FH}. The theorem is proved in Subsection~\ref{pf:FH}
by reducing it to a special case of finite-horizon dynamic programming
for controlled Markov processes \cite{flemingcontrolled2006}.  

\gap

\begin{theorem}{\it
    \label{thm:FHVerification}
    Assume that there is a function $V\in\D$ that satisfies:
    \begin{align}\label{eq:bellman}
      \inf_u 
      \left[
        c(s,x,u) + \A^u V(s,x)
      \right]&=0, 
\\
      \label{eq:boundary}
      V(S,x) &= \Psi(x), 
    \end{align}
where \eqref{eq:bellman} holds for all $(s,x,u)\in [0,S)\times \X\times \U$ and \eqref{eq:boundary} holds for all $x\in\X$.

  Then $V(s,x)\le J(s,x;U)$ for every feasible policy, $U$. 

    Furthermore, if a policy $U^*_r$ and associated state process $X_r^*$,
    with $X_s^*=x$, satisfy 
    \begin{equation*}
    U_r^* \in \argmin_u \left[
      c(r,X_r^*,u) + \A^u V(r,X_r^*)
    \right],
    \end{equation*}
for almost all $(r,\omega)\in [s,S]\times \Omega$,
    then $V(s,x) = J(s,x;U^*)$. In other words, $U^*_s$ is optimal.  
  }
\end{theorem}

\gap

\begin{example}
  Consider the problem of maximizing $\E\left[X_S^{\eta}\right]$, with
  $\eta \in (0,1)$ subject to the time-changed dynamics
  \begin{align*}
    dY_t & = U(\zeta_t) Y_t (\mu_1 dt + \sigma_1 dW_1(t)) \\
    & \quad\quad + (1-U(\zeta_t)) Y_t (\mu_2 dt
    + \sigma_2 dW_2(t)) \\ 
    X_s &= Y(\tau_s),
  \end{align*}
  where $W_1(t)$ and $W_2(t)$ are independent Brownian motions.  The
  problem can be interpreted as allocating wealth between stocks
  modeled by time-changed geometric Brownian motions: $Z_i(s) =
  R_i(\tau_s)$, where $dR_i(t) = R_i(t)(\mu_i dt + \sigma_i
  dW_i(t))$. 

  Let $u^*$ be the optimal solution and $\rho^*$ be the
  optimal value of the following quadratic maximization problem:
  \begin{multline*}
  \max_u \bigg[
    \frac{1}{2}\eta(\eta-1) \left(
      (u\sigma_1)^2 + ((1-u)\sigma_2)^2
    \right) \\
    + \eta \left(u \mu_1 + (1-u) \mu_2\right)
  \bigg].
  \end{multline*}
    If $\rho^*\in\dom(\beta)$, it can be verified by elementary stochastic calculus that
    $V(s,x)$ given by 
    \begin{equation*}
    V(s,x) = e^{\beta(\rho^*)(S-s)}x^{\eta}
    \end{equation*}
    satisfies the dynamic programming equations, \eqref{eq:bellman} and
    \eqref{eq:boundary}, with $\X \times \U = [0,\infty) \times \R$ and $\max$ replacing $\min$. The corresponding
    optimal input is $U_s^* = u^*$.
\end{example}

\subsection{Linear Quadratic Regulators}

In this section, Theorem~\ref{thm:FHVerification} is specialized to
linear systems with quadratic cost. The result (with no Brownian
forcing) was originally presented in \cite{lamperskitimechanged2013}, using a
proof technique specialized for linear systems.

\gap

\begin{problem}\label{prob:LQR}
  Consider linear dynamics 
  \begin{equation}\label{linSys}
    dY_t = (AY_t +BU(\zeta_t)) dt + MdW_t, 
  \end{equation}
  subject to the time change $X_s = Y(\tau_s)$. Here $\X = \R^n$ and
  $\U = \R^p$.

  The time-changed linear quadratic regulator problem over  time horizon $[0,S]$ is to find a policy $U_s$ that solves 
  \begin{equation*}
  \min_U \E\left[\int_0^S \left(X_s^\tp Q X_s + U_s^\tp R U_s\right) ds + X_S^\tp \Phi X_S \right],
  \end{equation*}
  over all c\`agl\`ad, $\F^{\zeta,X}$-adapted policies. Here $Q$ and $\Phi$ are positive semidefinite, while $R$ is positive definite. 
\end{problem}

\gap

The following lemma introduces the mappings used to construct the
optimal solution for the time-changed linear quadratic regulator
problem. The lemma is proved in Appendix~\ref{pf:mappings} by showing
that each mapping may be computed from $\beta(\tilde{A})$ for an
appropriately defined matrix $\tilde{A}$.

\gap

\begin{lemma}\label{lem:mappings} {\it
  Let $P$ be an $n\times n$ matrix. If $\{0\}\cup\spec(2A) \subset
  \dom(\beta)$, then the following linear mappings are well defined:
  \begin{align*}
    F(P) &= b(A^\tp P + PA) + \int_0^{\infty} \left(e^{A^\tp
        t}Pe^{At}-P\right)\lambda(dt) \\
    G(P) &= bP + \int_0^{\infty} e^{A^\tp t} P \int_0^t e^{Ar}dr
    \lambda(dt) \\
    H(P) &= \int_0^{\infty} \int_0^t e^{A^\tp r} dr P \int_0^t e^{A\rho}
    d\rho \lambda(dt) \\
    g(P) &= \Tr\left(
      P\left(
        b MM^\tp + \int_0^\infty\int_0^t e^{Ar} MM^\tp e^{A^\tp r} dr \lambda(dt) 
      \right)
    \right).
  \end{align*}
  Furthermore, $F$, $G$, and $H$ satisfy
  \begin{align*}
    \E\left[e^{A^\tp \tau_s} P e^{A \tau_s} \right] & = P + sF(P) + O(s^2) \\
    \E\left[e^{A^\tp \tau_s} P \int_0^{\tau_s} e^{Ar} dr\right] &= sG(P) + O(s^2)
    \\
    \E\left[
      \int_0^{\tau_s} e^{A^\tp r}dr P \int_0^{\tau_s} e^{A\rho} d\rho
    \right] &= s H(P) + O(s^2).
  \end{align*}
}
\end{lemma}

\gap

\begin{remark}{\rm
    The descriptions of $F$, $G$, and $H$ in terms of expectations are not required for the proof below. They are given to demonstrate that the formulas in terms of $(b,\lambda)$ coincide with the formulas from \cite{lamperskitimechanged2013}.
  }
\end{remark}

\gap

\begin{example}{\rm
    With no temporal noise, the mappings reduce to 
    \begin{equation}\label{eq:noiselessMat}
      \begin{array}{rclrcl}
      F(P) = A^\tp P+YP, & G(P)  = P, \\
      H(P)  = 0, & g(P) = \Tr(PMM^\tp ).
      \end{array}
    \end{equation}
    Furthermore, since $\beta(z)=z$ is analytic everywhere, these
    formulas are true for any state matrix, $A$. 
  }
\end{example}

\gap

\begin{example}{\rm
    Consider an arbitrary strictly increasing subordinator with Laplace
    exponent $\psi$. 
    Let  $A=\mu$ where $\mu$ is a real, non-zero scalar with
    $2\mu \in\dom(\beta)$. Let $M$ be a scalar. Combining \eqref{eq:laplace} with the formula $\int_0^te^{\mu \sigma}d\sigma =
    \mu^{-1}(e^{\mu t}-1)$ shows that
    \begin{align*}
      F(P) &= \beta(2\mu)P \\
      G(P) &= \mu^{-1}(\beta(2\mu)-\beta(\mu))P \\
      H(P) &= \mu^{-2}(\beta(2\mu)-2\beta(\mu))P \\
      g(P) &= \frac{1}{2}\mu^{-1} \beta(2\mu) M^2 P.
    \end{align*}
  }
\end{example}

\begin{theorem}\label{thm:LQR} {\it
  Say that $\{0\}\cup \spec(2A) \subset \dom(\beta)$. Define the function
  $V(s,x)=x^\tp P_s x + h_s$ by the backward differential equations 
  \begin{align*}
    -\frac{d}{ds} P_s &= Q + F(P_s) - G(P_s) B(R+B^\tp H(P_s)B)^{-1} B^\tp
    G(P_s)^\tp \\
    -\frac{d}{ds} h_s &= g(P_s),
  \end{align*}
  with final conditions $P_S = \Phi$ and $h_S = 0$. The function
  $V(s,x)$ satisfies  dynamic
  programming equations, \eqref{eq:bellman} and \eqref{eq:boundary}, and
  the optimal policy is given by 
  \begin{align*}
    U_s &= K_s X_{s^-} \\
    K_s &= -(R+B^\tp H(P_s)B)^{-1} B^\tp
    G(P_s)^\tp.
  \end{align*}
}
\end{theorem}

A straightforward variation on the proof of Theorem~\ref{thm:LQR} shows that for any linear policy, $U_s = L_s X_{s^-}$, the cost-to-go is given by
\begin{equation*}
J(s,x;U) = x^\tp Z_s x + p_s, 
\end{equation*}
where $Z_s$ and $k_s$ satisfy the backward differential equations 
\begin{align*}
  -\frac{d}{ds} Z_s &= Q+F(Z_s)+L_s^\tp B^\tp G(Z_s)^\tp +G(Z_s)BL_s 
  \\ &
  +L_s^\tp(R+B^\tp H(Z_s)B)L_s \\
  -\frac{d}{ds}p_s &= g(Z_s).
\end{align*}
In the following example, these formulas are used in order to compare the performance of the policy from Theorem~\ref{thm:LQR} with the policy $U_s = L_sX_{s^-}$, where $L_s$ is the standard LQR gain, not compensating for temporal noise. 

\gap

\begin{figure}
      \centering
      \includegraphics{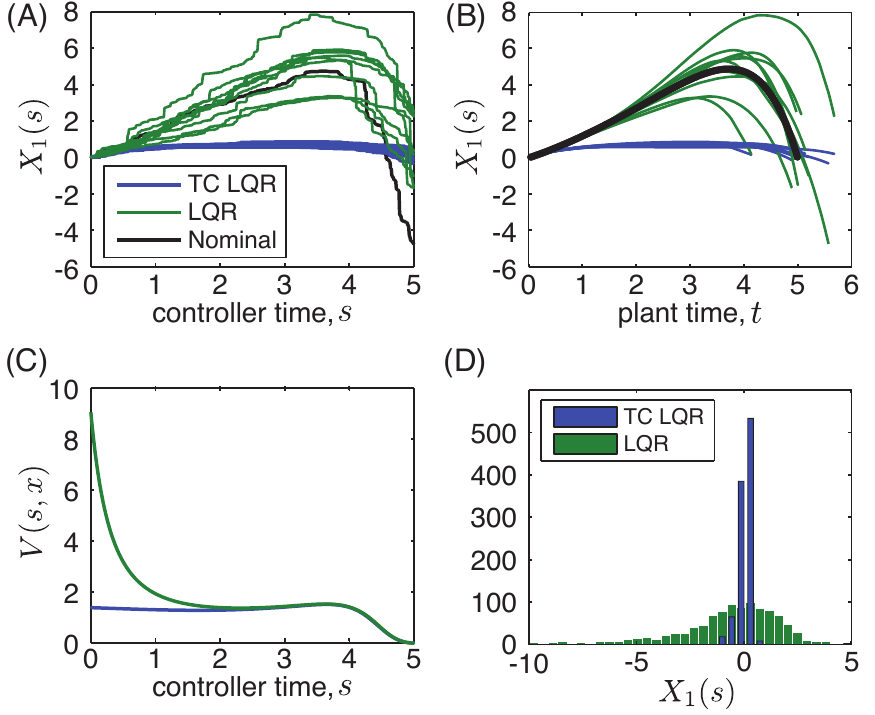}
      \caption{\label{LQROptCompare} (A) Plots of $X_1(s)$ under the
        optimal policy and the LQR policy for $10$ realizations of
        $\tau_s$. The initial condition is $x = [0,1]^\tp$. (B) The
        same plots under time variable $t$. The black line shows the
        LQR trajectory with no temporal noise. In the case of no
        temporal noise, the classical LQR uses high gains near $t=0$
        to produces high-speed trajectories such that $Y_1$ approaches
        $0$ at final time. In this case, timing errors lead to wide
        variation in the final position.  The optimal policy reduces
        the speed of the trajectory near $s=0$ order to minimize the
        effects of temporal noise.  (C) The optimal cost
        $V(s,x)$ and $J(s,x;U)$ for the LQR policy are plotted for
        $x=[0,1]^\tp$. As expected, $V(s,x)\le J(s,x;U)$. Furthermore,
        as the time-horizon increases, the LQR policy depends strongly
        on timing information, and so temporal noise leads to higher
        cost as $s$ goes to $0$. (D) A histogram of
        the final positions, $X_1(S)$, for $1000$ realizations of
        $\tau_s$. The optimal controller leads to $X_1(S)$ being
        tightly distributed around $0$, while the LQR controller gives
        a wide spread of $X_1(S)$ values. The errors in the final
        position lead to increased cost for the LQR controller.  }
    \end{figure}

\begin{example}{\rm
    Consider the system defined by the state matrices
    \begin{equation*}
    A=\begin{bmatrix}
      0.75 & 1 \\
      0 & 0.75
    \end{bmatrix},
    \qquad
    B=\begin{bmatrix}
      0 \\ 1
    \end{bmatrix}, \qquad M=0,
    \end{equation*}
    with cost matrices  given by 
    \begin{equation*}
    R = 0.5,\qquad Q=0, \qquad 
    \Phi = \begin{bmatrix}
      1 & 0 \\
      0& 0 
    \end{bmatrix}.
    \end{equation*}
    Let $\tau_s$ be the inverse Gaussian subordinator with $\gamma=\delta=2$. The condition, $\spec(2A)\subset\dom(\beta)$, is satisfied since $2\cdot 0.75 = 1.5 < \gamma^2/2 = 2$. Figure~\ref{LQROptCompare} compares the optimal policy with the standard LQR policy. 
  }
\end{example}

\section{Proofs of Main Results}\label{sec:pf}

\subsection{Proof of Theorem~\ref{timeChangedDiffusion}}
\label{pf:tcDiffusion}
From the definition of $X_s$,
\begin{equation}\label{XSolDef}
  X_s = X_0 + \int_0^{\tau_s}F_t dt + \int_0^{\tau_s}G_t dW_t. 
\end{equation}
Thus, provided that \eqref{timeChange} holds, claims \ref{finiteVariation})
and \ref{martingale}) imply that $X_s$ must be an $\F^{\tau,W}$
semimartingale. The claims are proved as follows. 
\begin{equation*}
\mathrm{Var}\left(
  \int_0^{\tau_s} F_t dt 
\right) \le \int_0^{\tau_s} \|F_t\|dt <\infty \textrm{ almost surely.} 
\end{equation*}
Therefore \ref{finiteVariation}) holds.

To prove \ref{martingale}), note that for $0\le r \le s$ we have
\begin{align*}
\E\left[
  \int_0^{\tau_s}G_t dW_t
  \:\vline\: 
  \F^{\tau,W}_r
\right] &=
\int_0^{\tau_r}G_t dW_t + \E\left[
  \int_{\tau_r}^{\tau_s}G_t dW_t
  \:\vline\: 
  \F^{\tau,W}_r
\right] 
\\
&= \int_0^{\tau_r}G_t dW_t. 
\end{align*}
Furthermore,
\begin{equation}\label{martingaleBound}
  \E\left[
    \left\|
      \int_0^{\tau_s}G_t dW_t
    \right\|
  \right]^2
  \le  \E\left[
    \left\|
      \int_0^{\tau_s}G_t dW_t
    \right\|^2
  \right]
  < \infty ,
\end{equation}
where the inequality follows from Jensen's inequality.
Thus \ref{martingale})
holds. 

Now \eqref{timeChange} must be proved. For more compact notation, define the processes $H_t$ and $Z_t$ as 
\begin{equation*}
H_t = \begin{bmatrix} F_t & G_t\end{bmatrix} \qquad 
Z_t = \begin{bmatrix}
  t \\ W_t 
\end{bmatrix}
\end{equation*}
so that $X_s$ may be written as 
\begin{equation*}
X_s = \int_0^{\tau_s} H_t dZ_t.
\end{equation*}

Note that $Z(\tau_s) = [\tau_s,W(\tau_s)^\tp ]^\tp$. Since $\tau_s$ is
a subordinator, $W(\tau_s)$ is a L\'evy process on $\R^d$, with L\'evy
symbol
\begin{equation*}
\eta_{W_{\tau}}(z) = -\frac{1}{2}b z^\tp z +
\int_{\R^d}(e^{iz^{\tp} x}-1)\mu_{W,\tau}(dt),
\end{equation*}
for some L\'evy measure $\mu_{W,\tau}$. (See Theorem 1.3.25 and
Theorem 1.3.33, respectively, in \cite{applebaumlevy2004}.) Thus, the
continuous part of $W(\tau_s)$ is a Brownian motion with
$\E\left[W(\tau_s)W(\tau_s)^\tp\right] = b I$. 

Define the $\tilde{Z}_s$  by removing the jumps from $Z(\tau_s)$.
\begin{equation*}
\tilde{Z}_s = Z(\tau_s) - \sum_{0\le r \le s} \left(Z(\tau_r)-Z(\tau_r^-)\right).
\end{equation*}
It follows that $\tilde{Z}_s = [bs,\sqrt{b} \tilde{W}_s^\tp]^\tp$,
where $\tilde{W}_s$ is the Brownian motion from the theorem statement.
Thus, \eqref{timeChange} can be equivalently written as 
\begin{equation}\label{compactIntegral}
  X_s = \int_0^s H(\tau_{r^-})d\tilde{Z}_s + \sum_{0\le r \le s}
  \int_{\tau_{r^-}}^{\tau_r} H_s dZ_t
\end{equation}

Now \eqref{compactIntegral} will be evaluated. If $b=0$, then
$\tilde{Z}_s = 0$ and  $\tau_s
= \sum_{0\le r \le s} \tau_{r} - \tau_{r^-}$. Thus, 
\begin{equation*}
X_s = \sum_{0\le r \le s} \int_{\tau_{r^-}}^{\tau_r} H_s dZ_t,
\end{equation*}
so, in this case, \eqref{timeChange} holds.

Now assume $b>0$. The cases when $\tau_s$ has finite rate
($\lambda((0,\infty)) < \infty$) and infinite rate
($\lambda((0,\infty)) = \infty$) will be treated separately.

\textbf{Finite Rate.} Let $r_0=0$ and let $r_1,r_2,\ldots$ be the jump times of
$\tau_s$. With probability $1$, there exist a finite (random) integer
$L$ such that $L$ jumps occur over $[0,s]$.  
Note that \eqref{compactIntegral} may be expanded as 
\begin{align}\label{finiteRateIntermediate}
  X_s &= \int_{\tau(r_L)}^{\tau_s} H_t dZ_t 
\\ & \nonumber
+ \sum_{k=0}^{L-1} 
  \left[
    \int_{\tau(r_k)}^{\tau(r_{k+1}^-)} H_t dZ_t +
    \int_{\tau(r_{k+1}^-)}^{\tau(r_{k+1})} H_t dZ_t 
  \right]
\end{align}
Let $s_0^n \le s_1^n \le \cdots \le  s_{K_n}^n$ be a sequence of partitions such that 
\begin{equation*}
\begin{array}{lc}
  \lim_{n\to\infty} s_{K_n}^n = \infty & a.s. \\
  \lim_{n\to\infty}\sup\{|s_{k+1}^n - s_k^n| : k=0,\ldots ,K_{n}-1\} = 0 & a.s. \\
  \{r_i: r_i \le s_{K_n}^n\} \subset \{s_0^n,\ldots,s_{K_n}^n\}.
\end{array}
\end{equation*}
The last condition ensures that the jump times are contained in the
partition. 

Note that between jumps (i.e. $s\in
[r_k,r_{k+1})$), $\tau_s = bs + \tau^d(r_k)$, where $\tau^d_s$ is the
discontinuous part of $\tau_s$. Since $b>0$ follows that the sequence
$\tau(s_0^n),\tau(s_1^n),\ldots$, satisfies the following properties, almost surely:
\begin{equation*}
\begin{array}{l}
  \lim_{n\to\infty} \tau(s_{K_n}^n) = \infty 
\\
  \lim_{n\to\infty} \sup\{
  |\tau(s_{i+1}^n) - \tau(s_i^n)|: \exists k \:\: s.t.\:\: r_k \le s_i^n
  < r_{k+1} \} = 0 
\end{array}
\end{equation*}
Using a standard argument from stochastic integration (see Theorem II.21 of
\cite{protterstochastic2004}), the integral from $\tau(r_k)$ to
$\tau(r_{k+1}^-)$ may be evaluated as 
\begin{align}
  \nonumber 
\MoveEqLeft
  \int_{\tau(r_k)}^{\tau(r_{k+1}^-)} H_s dZ_t 
\\ \nonumber
&=
  \lim_{n\to\infty} \sum_{r_k \le s_i^n < r_{k+1}} H(\tau(s_i^n)) 
  \left(Z(\tau(s_{i+1})) - Z(\tau(s_i))\right) \\
  \nonumber
  &= \lim_{n\to\infty} \sum_{r_k \le s_i^n < r_{k+1}} H(\tau(s_i^n)) 
  \left(\tilde{Z}(s_{i+1}) - \tilde{Z}(s_i)\right) \\
  & = \int_{r_k}^{r_{k+1}} H(\tau_{s^-}) d\tilde{Z}_s.
  \label{finiteRateApprox}
\end{align}
The second equality uses the fact that no jumps occur over
$(r_k,r_{k+1})$. The result now follows by combining
\eqref{finiteRateIntermediate} and \eqref{finiteRateApprox}.

\textbf{Infinite Rate.} Let $\epsilon_n>0$ be a sequence decreasing to
$0$, at a rate to be specified later. Define $\tau_s^n$ to be the
process by removing all jumps of size at most $\epsilon_n$ from
$\tau_s$:
\begin{equation}\label{tauNdef}
  \tau_s^n = bs + \int_{\epsilon_n}^\infty t N(s,dt).
\end{equation}
Let $r_0^n=0$, and let $r_1^n,r_2^n,\ldots$ be the jump times of
$\tau^n_s$. Let $L_s^n = \sup\{k:r_k^n \le s\}$. With probability $1$,
$L_s^n < \infty $. If $\epsilon_n$ are chosen as in Lemma~\ref{lem:smallTimes} from
Appendix~\ref{sec:lemmas}, then $X_s$ may be computed as a limit
\begin{align}\label{infiniteRateFirstExpand}
\MoveEqLeft
  X_s = \lim_{n\to\infty}\left[\int_{\tau(r_{L_s^n}^n)}^{\tau_s}H_t dZ_t  \right.
  \\ \nonumber &
    +\sum_{k=0}^{L_s^n-1}
      H(\tau(r_k^n))
      \left(
        Z(\tau(r_{k+1}^{n-})) - Z(\tau(r_k^n))
      \right)
\\ \nonumber &
\left.
      +
\sum_{k=0}^{L_s^n-1}
\int_{\tau(r_{k+1}^{n-})}^{\tau(r_{k+1}^n)} H_t dZ_t
  \right].
\end{align}
Note that $Z(\tau(r_{k+1}^{n-})) - Z(\tau(r_k^n))$ may be expressed as 
\begin{align*}
\MoveEqLeft
Z(\tau(r_{k+1}^{n-})) - Z(\tau(r_k^n)) 
\\ &=
\tilde{Z}(r_{k+1}^n)-\tilde{Z}(r_k^n)
+ \sum_{\substack{r_k^n < r \le r_{k+1}^n \\ \Delta \tau_r \le
    \epsilon_n}}
\left(
  Z(\tau_r) - Z(\tau_{r^-})
\right).
\end{align*}

Note that the terms in the summation all vanish as $\epsilon_n \to
0$. Furthermore, $r_{L_s^n}^n \uparrow s$, almost surely.  Thus,
\eqref{infiniteRateFirstExpand} can be expressed as
\begin{align*}
  X_s &= \lim_{n\to\infty}\left[\sum_{k=0}^{L_s^n-1}
      H(\tau(r_k^n))
\left(
        \tilde{Z}(r_{k+1})-\tilde{Z}(r_k)
\right)
      \right.
\\ & \left. 
      +
\sum_{k=0}^{L_s^n-1}
\int_{\tau(r_{k+1}^{n-})}^{\tau(r_{k+1}^n)} H_t dZ_t
  \right],
\end{align*}
and the result now follows using Theorem II.21 of
\cite{protterstochastic2004}. \hfill\IEEEQED

\subsection{Proof of Theorem \ref{thm:FHVerification}}
\label{pf:FH}

Theorem~\ref{thm:FHVerification} is a special case of finite-horizon
dynamic programming for controlled Markov processes (Theorem
III.8.1 of \cite{flemingcontrolled2006}), provided that the following
two conditions hold for all $h\in \D$:
\begin{enumerate}[(i)]
\item The backward evolution operator, defined in \eqref{generatorLim}
  is given by the formula in \eqref{generator}.
\item If 
  \begin{gather*}
    \E\left[|h(S,X_S)| \: | \: X_s=x\right] < \infty
    \intertext{and}
    \E\left[\int_s^S \left|\A^{U_r} h(r,X_r)\right|dr \:\vline \:
    X_s=x\right] <\infty,
  \end{gather*}
  then the Dynkin formula holds:
  \begin{multline}
    \label{Dynkin}
      \E\left[h(S,X_{S})\:\vline\:
        X_s=x\right] - h(s,x) 
      \\ 
      = 
      \E\left[
        \int_s^{S} \A^{U_r}h(r,X_r)dr \:\vline\:
        X_s = x
      \right].
  \end{multline}
\end{enumerate}

First, using Theorem~\ref{timeChangedDiffusion}, a more explicit
formula for $X_s$ is derived, and then using It\^o's formula for
semimartingales, a formula for $h(s,X_s)$ is given. Using the formula
for $h(s,X_s)$, equations \eqref{generatorLim} and \eqref{Dynkin} are
then proved.

Note that $Y(\tau_{r^-})=X_{r^-}$ and for all $t\in
[\tau_{r^-},\tau_r]$, $\zeta_t = r$. Therefore 
\begin{gather*}
  F(\zeta(\tau_{r^-}),Y_{\tau_{r^-}},U(\zeta(\tau_{r^-}))) = 
  F(r,X_{r^-},U_r) \\
  F(\zeta_t,Y_t,U(\zeta_t)) = 
  F(r,Y_t,U_r) \qquad \textrm{for all} \qquad t\in [\tau_{r^-},\tau_r].
\end{gather*}
The expressions for $G$ are similar. Thus,
Theorem~\ref{timeChangedDiffusion} implies that $X_s$ is given by
\begin{multline}
    \label{TCDynamics}
  X_s = X_0 + b\int_0^s F(r,X_{r^-},U_{r}) dr + \sqrt{b}\int_0^s G(r,X_{r^-},U_{r})d\tilde{W}_r 
  \\
  + \sum_{0\le r \le s}
  \left(
    \int_{\tau_{r^-}}^{\tau_r}F(r,Y_t,U_{r})dt + 
    \int_{\tau_{r^-}}^{\tau_r}G(r,Y_t,U_{r})dW_t
  \right).
\end{multline}

Now a formula for $h(s,X_s)$ will be derived. Note that for any c\`agl\`ad,
$\F^{\tau,W}$-adapted process, $Z_s$,   the stochastic integral with respect to $X_s$ is given by 
\begin{multline*}
  \int_0^s Z_r dX_r = \int_0^s Z_r b F(r,X_{r^-},U_r)dr
  \\
  +\int_0^s Z_r \sqrt{b} G(r,X_{r^-},U_r) dW_r + 
  \sum_{0\le r\le s} Z_r (X_{r}-X_{r^-}). 
\end{multline*}
Furthermore, the continuous part of the quadratic variation is given by
\begin{equation*}
[X,X]_s^c = \int_0^s\frac{1}{2}b G(r,X_{r^-},U_r)G(r,X_{r^-},U_r)^\tp dr. 
\end{equation*}
Thus It\^o's formula for semimartingales (see \cite{protterstochastic2004}) implies that $h(s,X_s)$ is given by
\begin{multline}\label{TCIto}
  h(s,X_s)= h(0,X_0) 
  \\
  + \int_0^s \left(
    \frac{\partial h(r,X_{r^-})}{\partial r}
    +b\frac{\partial h(r,X_{r^-})}{\partial x} F(r,X_{r^-},U_{r}) 
  \right) dr \\
  +\int_0^s\frac{1}{2}b\Tr\left(
    G(r,X_{r^-},U_{r})^\tp 
    \frac{\partial^2 h(r,X_{r^-})}{\partial x^2} G(r,X_{r^-},U_{r})
  \right) dr 
  \\
  +\int_0^s \sqrt{b}\frac{\partial h(r,X_{r^-})}{\partial x}
  G(r,X_{r^-},U_{r}) d\tilde{W}_r 
  \\
  +
  \sum_{0\le r \le s}
  \left(
    h(r,X_r) - h(r,X_{r^-})
  \right).
\end{multline}

Now \eqref{generatorLim} will be derived from \eqref{TCIto}. Assume that $X_s=x$ and $U_r = u$ for $r\in [s,s+h]$. Then \eqref{TCIto} implies that 
\begin{multline}
  \label{SmallItoSum}
  \E[h(s+h,X_{s+h})]=  h(s,x)
  \\
  +
  \E\left[
    \int_s^{s+h} 
    \left(
      \frac{\partial h(r,X_{r^-})}{\partial r}
      +b\frac{\partial h(r,X_{r^-})}{\partial x} F(r,X_{r^-},u) \right)dr
  \right]
  \\
  +
  \E\left[
    \frac{1}{2}b\Tr\left(
      G(r,X_{r^-},u)^\tp 
      \frac{\partial^2 h(r,X_{r^-})}{\partial x^2} G(r,X_{r^-},u)
    \right) dr
  \right]
  \\
  +
  \E\left[
    \sum_{s< r \le s+h}
    \left(
      h(r,X_r) - h(r,X_{r^-})
    \right)
  \right].
\end{multline}

If $r>s$, the Brownian motions $W_t$ for $t\in [\tau_{r^-},\tau_r]$
and $\hat{W}_t$ for $t\in [0,\Delta\tau_r]$ are identically
distributed and independent of $\F_s^{\tau,W}$. Therefore, using
\eqref{constUDrift}, and given that $X_s=x$ and $U_r=u$, the
expectations of the jump terms may be written as
\begin{align*}
\MoveEqLeft
\E\left[
  h(r,X_r)-h(r,X_{r^-})
\right] =
\\ &
\E\left[
  \E_{\hat{W}}\left[h\left(r,Y_{r,\Delta\tau_r}^{X_{r^-}u}\right)\right]-h(r,X_{r^-})
\right].
\end{align*}
Thus, the term at the bottom of \eqref{SmallItoSum} may be evaluated as a Poisson integral:
\begin{align*}
  \MoveEqLeft[0]
  \E\left[
    \sum_{s< r \le s+h}
    \left(
      h(r,X_r) - h(r,X_{r^-})
    \right)
  \right] 
  \\
  &=
  \E\left[\int_s^{s+h}\int_0^{\infty}
    \left(
      \E_{\hat{W}}\left[h\left(r,Y_{rt}^{X_{r^-}u}\right)\right]-h(r,X_{r^-})
    \right)N(dr,dt)
  \right]
  \\
  &=\E_W\left[
    \int_s^{s+h} \int_0^{\infty}\left(
      \E_{\hat{W}}\left[h\left(r,Y_{rt}^{X_{r^-}u}\right)\right]-h(r,X_{r^-})
    \right)\lambda(dt)dr
  \right]
\end{align*}
where the second is equation is justified by Fubini's theorem and \eqref{bddCond}. 

By evaluating the limit on the right side of 
\eqref{generatorLim}, the formula in \eqref{generator} is recovered.     

Turning to \eqref{Dynkin}, since $X_s$ is $\F_{s}^{\tau,W}$
measurable, it suffices to prove that
\begin{equation*}
\E\left[
  h(S,X_S)-h(s,X_s)-\int_s^S \A^{U_s} h(r,X_r)dr 
  \:\vline\:
  \F_s^{\tau,W}
\right]=0.
\end{equation*}
Since $X_r(\omega)=X_{r^-}(\omega)$ for almost all $(r,\omega)\in
[s,S]\times \Omega$, it follows that 
\begin{align*}
\MoveEqLeft
\E\left[
  \int_s^S \A^{U_s} h(r,X_r)dr 
  \:\vline\:
  \F_s^{\tau,W}
\right] = 
\\ &
\E\left[
  \int_s^S \A^{U_s} h(r,X_{r^-})dr 
  \:\vline\:
  \F_s^{\tau,W}
\right]. 
\end{align*}

Combining \eqref{generator}
and \eqref{TCIto} and, for brevity, omitting $\F_s^{\tau,W}$ from the
expectations, implies that
\begin{multline}
  \label{JumpExp}
  \E\left[
    h(S,X_{S}) - h(s,X_s) - \int_s^{S} \A^{U_{r}} h(r,X_{r^-}) dr  \right] = 
  \\
  \hspace{-5pt}
  \E\left[
    \sum_{s<r\le S} \left(
      h(r,X_r)-h(r,X_{r^-})
    \right) \right]-\\
  \E\left[\int_s^{S}
    \int_0^{\infty}
    \left(
      \E_{\hat{W}}[h(r,Y_t^{X_{r^-}U_{r}})]-h(r,X_{r^-})
    \right)\lambda(dt)
    dr
  \right],
\end{multline}
As in the proof of \eqref{generatorLim}, the two terms at the bottom
of \eqref{JumpExp} are equal in expectation. Thus \eqref{Dynkin} holds
and the proof is complete.  \hfill\IEEEQED

\subsection{Proof of Theorem~\ref{thm:LQR}}
Assume $V(s,x)= x^\tp P_s x + h_s$. Applying the backward evolution
operator corresponding to \eqref{linSys} to $V(s,x)$ results in
\begin{align}
  \label{linGenerator}
  \MoveEqLeft[0]
  \A^u V(s,x) = x^\tp \dot{P}_s x + \dot{h}_s + 
  bx^\tp (A^\tp P_s + P_s A) x + 2b x^\tp P_sBu 
  \\ \nonumber  
  &+
b\Tr\left(P_s MM^\tp\right)+
\int_0^\infty\left(\E_{\hat{W}}\left[{Y_t^{xu}}^\tp P_s
      Y_t^{xu}\right] 
    -x^\tp P_s x\right)\lambda(dt).
\end{align}

Note that $\E_{\hat{W}}\left[Y_t^{xu} {Y_t^{xu}}^\tp\right]$ may be
evaluated at 
\begin{equation*}
\E_{\hat{W}}\left[Y_t^{xu} {Y_t^{xu}}^\tp\right] =
y_t^{xu}{y_t^{xu}}^\tp + \Sigma_t,
\end{equation*}
where $y_t^{xu}$ is the mean of $Y_t^{xu}$ and $\Sigma_t$ is the
covariance. A standard argument in linear stochastic differential
equations shows that the mean and covariance are given by
\begin{align*}
  y_t^{xu} &= e^{At} x + \int_0^t e^{Ar}dr Bu, \\
  \Sigma_t &= \int_0^t e^{Ar} MM^\tp e^{A^\tp r} dr.
\end{align*}
Thus, the integral in \eqref{linGenerator} may be written as 
\begin{align}
  \label{integralExplicit}
  \MoveEqLeft
  \int_0^\infty\left(\E_{\hat{W}}\left[{Y_t^{xu}}^\tp P_s
      Y_t^{xu}\right] 
    -x^\tp P_s x\right)\lambda(dt)=
  \\ & \nonumber
  \int_0^{\infty}
  \begin{bmatrix}
    x \\ u
  \end{bmatrix}^\tp 
  \begin{bmatrix}
    \hat F(t,P_s) & \hat G(t,P_s) B \\
    B^\tp \hat G(t,P_s)^\tp & B^\tp \hat H(t,P_s) B
  \end{bmatrix}
  \begin{bmatrix}x \\ u\end{bmatrix}
  \lambda(dt)
  \\ \nonumber & 
  +
  \int_0^{\infty}
  \int_0^t\Tr\left(P_s e^{Ar}MM^\tp e^{A^\tp r}\right)dr
  \lambda(dt),
\end{align}
where the matrices $\hat F$, $\hat G$, and $\hat H$ are defined by
\begin{align*}
\hat F(t,P) &= e^{A^\tp t} Pe^{At} - P, \\\
\hat G(t,P) &=  e^{A^\tp t} P \int_0^t e^{A\rho}d\rho, \\
\hat H(t,P) &= \int_0^t e^{A^\tp r}dr P
    \int_0^t 
    e^{A\rho}d\rho.
\end{align*}
Combining \eqref{linGenerator} and \eqref{integralExplicit}, and using
the linear operators from Lemma~\ref{lem:mappings} gives
\begin{equation*}
\A^u V(s,x) = \begin{bmatrix} x\\ u\end{bmatrix}^\tp 
\begin{bmatrix}
  \dot{P}_s + F(P_s) & G(P_s) B \\
  B^\tp G(P_s) & B^\tp H(P_s) B 
\end{bmatrix}
\begin{bmatrix} x\\ u\end{bmatrix} + \dot{h}_s + g(P_s).
\end{equation*}
Therefore, adding the cost gives
\begin{align*}
  \MoveEqLeft
  x^\tp Qx +u^\tp R u + \A^u V(s,x) = \dot{h}_s + g(P_s)
  \\&
  +\begin{bmatrix} x\\ u\end{bmatrix}^\tp 
  \begin{bmatrix}
    Q+\dot{P}_s + F(P_s) & G(P_s) B \\
    B^\tp G(P_s) & R + B^\tp H(P_s) B 
  \end{bmatrix}
  \begin{bmatrix} x\\ u\end{bmatrix}.
\end{align*}
The result now follows from quadratic minimization.
\hfill\IEEEQED

\section{Discussion}\label{sec:discussion}
The work in this paper lays a theoretical foundation for future
research on biological motor control
\cite{todorovoptimal2002,harrissignal-dependent1998}, finance
\cite{rogersoptimal2013}, and multi-agent control \cite{jadbabaiecoordination2003,olfatisaberconsensus2004} in context that the
controller is uncertain about the time of the plant.  To reason about
these problems, theoretical extensions will include
time estimation from sensory data, optimal control control with
different time horizons, and control with multiple noisy clocks.

\subsection*{Time Estimation}
This paper focuses on state feedback problems, and the optimal
solution from dynamic programming only depends on the value of the
current state. For other problems, an explicit estimate of time could
be valuable. In option pricing, inferences about the ``business time''
can be used to estimate the volatility of stock prices
\cite{gemanstochastic2002}. To perceive time, humans appear to
integrate sensory cues about the passage of time in a Bayesian manner
\cite{shibayesian2013}; humans also appear to incorporate sensory
information about the timing of events to improve state estimation
\cite{ankaralihaptic2014}.

In order to obtain a well-posed estimation problem, plant time
($\tau_s$) can be modeled as an unmeasured state. Particle filtering
techniques to estimate $\tau_s$ are currently being developed. 

In dance and music performance, movements are coordinated with an
estimate of time. The work in 
this paper will be combined with filtering methods to model
experimental data on timed movements. 

\subsection*{Variations on the Horizon}
This paper studies a controller horizon $[0,S]$, which is an interval
of time with respect to the clock measured by the controller. For
portfolio optimization, in which the controller measures calendar
time, such a horizon is sensible.  In human movements, however,
different tasks call for different time horizons. When reaching to an
object, a natural horizon would be the stopping time describing when
the object is touched. For rhythmic movements coordinated with external stimuli such as a metronome, 
a horizon over real time might
be sensible.

\subsection*{Multiple Clocks}

In this paper, it is assumed that the plant dynamics evolve according
to one clock, while the controller can measure a different clock. If
the plant consists of numerous subsystems, then each could potentially
evolve according to a different clock. This scenario arises in
portfolio problems, in which the goal is to allocate wealth between a
bank process which accrues interest at known, fixed rate and a stock process that evolves in a variable rate market
\cite{cvitanicoptimal2008}. Here, the bank process may be interpreted as evolving with respect to a perfect clock, while the stock process may be viewed as evolving with respect to a noisy clock.  In
engineering applications, such as mobile sensor networks, multiple
autonomous agents with their own clocks solve cooperative control
problems. Currently, problems arising from drift between clocks are
mitigated by using expensive clocks and time synchronization
protocols. The work in this paper will be extended to reduce the need
for precision timing and synchronization.

\section{Conclusion}\label{sec:conclusion}
This paper gives basic results on control with uncertainty in
time. The technical backbone of the paper is
Theorem~\ref{timeChangedDiffusion} which expresses the original plant
dynamics in terms of the controller's clock index. 
Using the new representation, the system becomes a controlled Markov
process, and thus existing dynamic programming theory can be applied. 
Given the dynamic programming equations, time changed versions
of linear quadratic control and a nonlinear portfolio problem problem are
solved explicitly. 


\appendices

\section{A Technical Lemma for Theorem~\ref{timeChangedDiffusion}}\label{sec:lemmas}

\begin{lemma}\label{lem:smallTimes} {\it
  Let $\tau_s$ be an infinite rate subordinator. Let $r_0^n=0$ and let $r_1^n\le
  r_2^n\le \cdots$
  be the jump times of $\tau_s^n$, from \eqref{tauNdef}. There is a
  sequence $S_n\to \infty$ and a sequence $\epsilon_n\downarrow 0$
  such that the following limits hold, almost surely
  \begin{align}
    \label{finalRLim}
    \lim_{n\to\infty}\sup\{r_i^n:r_i^n\le S_n\} &= \infty 
    \\
    \label{finalTauLim}
    \lim_{n\to\infty} \sup\{\tau_{r_i^n}:r_i^n \le S_n\}& = \infty
    \\
    \label{interJumpDomLim}
    \lim_{n\to\infty} \sup\{r_{i+1}^n-r_i^n : r_i^n \le S_n\} &= 0 
    \\
    \label{interJumpRanLim}
    \lim_{n\to\infty}\sup\{\tau_{r_{i+1}^{n-}}-\tau_{r_i^n}: r_i^n \le S_n\} & =0. 
  \end{align}
}
\end{lemma}
\begin{IEEEproof}
  First it will be shown that for any sequence $S_n>0$, a sequence
  $\epsilon_n\downarrow 0$ can be chosen such that
  \eqref{interJumpDomLim} and \eqref{interJumpRanLim} hold.  Then it
  will be shown how to choose $S_n$ so that \eqref{finalRLim} and
  \eqref{finalTauLim} hold. 

  Consider \eqref{interJumpDomLim}. Using Borel's lemma, it is
  sufficient to prove that for some constant $R>0$, and $\epsilon_n$
  sufficiently small,
  \begin{equation}
    \label{interJumpDomBound}
    \P\left(
      \sup_{r_i^n\le S_n}
      |r_{i+1}^n-r_i^n| \ge \frac{1}{2^n}
    \right) < \frac{R}{2^n}.
  \end{equation}

  For ease of notation, the superscripts on $r_i^n$ and the subscripts on
  $\epsilon_n$ and $S_n$ will be dropped. 

  With probability $1$, $\tau_s^n$ has only a finite number of jumps
  over $[0,S]$, so let $K=\max\{i:r_i\le S\}$. 

  Consider \eqref{interJumpDomBound}. Define the function $g(\epsilon)$ by
  \begin{equation*}
  g(\epsilon) = \int_{\epsilon}^\infty\lambda(dt).
  \end{equation*} 
  Note that the differences are $r_{i+1}-r_i$ are exponential random
  variables with rate parameter $g(\epsilon)$. Thus, the event
  that $r_{i+1}-r_i\ge 1/2^n$ is a Bernoulli random variable with
  probability $p(\epsilon)$ given by
  \begin{equation*}
  p(\epsilon)=\P\left(r_{i+1}-r_i\ge \frac{1}{2^n}\right) =
  g(\epsilon)\int_{1/2^n}^{\infty}e^{-g(\epsilon)x}dx = 
  e^{-g(\epsilon)/2^n}.
  \end{equation*}
  Let $J$ be the geometric random variable defined by
  \begin{equation*}
  J = \min\left\{i:r_{i+1}-r_i \ge \frac{1}{2^n}\right\}
  \end{equation*}
  Then the probability of $J$ is given by
  \begin{equation*}
  \P(J=k) = (1-p(\epsilon))^kp(\epsilon).
  \end{equation*}
  Using the definition of $K$ and $J$, the probability in
  \eqref{interJumpDomBound} may be written as 
  \begin{equation*}
  \P\left(
    \sup_{r_i\le S}
    |r_{i+1}-r_i| \ge \frac{1}{2^n}
  \right) = \P(J\le K).
  \end{equation*}
  Furthermore, given any constant $M>0$,
  \begin{equation}\label{DomBoundSplit}
    \P(J\le K) \le \P(J\le M) + \P(M\le K).
  \end{equation}
  Thus, \eqref{interJumpDomBound} may be bounded by bounding the terms
  on the right of \eqref{DomBoundSplit} separately. 

  Now $\P(M\le K)$ will be bounded. Note that $K$ is a Poisson random
  variable with parameter $Sg(\epsilon)$. Markov's inequality thus shows
  that 
  \begin{equation}\label{markovBound}
    \P(M\le K) \le \frac{1}{M} \E[K] = \frac{Sg(\epsilon)}{M}
  \end{equation}

  The term $\P(J\le M)$ can be computed exactly as 
  \begin{equation*}
  \P(J\le M) = p(\epsilon)\sum_{k=0}^M(1-p(\epsilon))^k = 1-(1-p(\epsilon))^{M+1}. 
  \end{equation*}
  Thus, \eqref{interJumpDomBound} will hold if $M$ can be chosen such
  that 
  \begin{equation}\label{firstMBounds}
    Sg(\epsilon)/M < 1/2^n
    \quad \textrm{and} \quad 1-(1-p(\epsilon))^{M+1} < 1/2^n. 
  \end{equation}
  Rearranging terms, \eqref{firstMBounds} is equivalent to 
  \begin{equation*}
  2^nS g(\epsilon) < M < \frac{\log\left(1-\frac{1}{2^n}\right)}{
    \log(1-p(\epsilon))}-1.
  \end{equation*}
  Therefore, a suitable constant $M$ exists if 
  \begin{equation}\label{eliminatedM}
    \left(2^nSg(\epsilon)+1\right)
    \log(1-p(\epsilon)) > \log\left(1-\frac{1}{2^n}\right).
  \end{equation}

  It will be shown that \eqref{eliminatedM} holds provided that
  $\epsilon$ is sufficiently small. Since $\tau_s$ has infinite rate,
  $\lim_{\epsilon\to 0} g(\epsilon) = \infty$. Thus, the limit of the
  left side of \eqref{eliminatedM} may be evaluated by L'H\^ospital's rule:
  \begin{align*}
    \lim_{\epsilon\to 0}\left(2^nSg(\epsilon)+1\right)
    \log(1-p(\epsilon)) &= \lim_{g\to\infty} 
    \frac{\log(1-e^{-g/2^n})}{\left(2^nSg+1\right)^{-1}} \\
    &=
    \lim_{g\to\infty}
    \frac{
      \frac{
        e^{-g/2^n}/2^n
      }{
        1-e^{-g/2^n}
      }
    }
    {
      -2^nS\left(2^nSg+1\right)^{-2}
    }
    \\
    &= -\frac{1}{4^nS} \lim_{g\to\infty}
    \frac{
      \left(2^nSg+1\right)^2
    }{
      e^{g/2^n}-1
    }
    \\
    &=0.
  \end{align*}
  Thus, when $\epsilon$ is sufficiently small, \eqref{interJumpDomBound}
  must hold.

  Now consider \eqref{interJumpRanLim}. Note that
  $\tau_{r_{i+1}^{n-}}-\tau_{r_i^n}$ can be expressed as
  \begin{equation*}
  \tau_{r_{i+1}^{n-}}-\tau_{r_i^n} = b(r_{i+1}^n-r_i^n)
  +\sum_{\substack{
      r_i^n < u \le r_{i+1}^n \\
      \Delta \tau_u \le \epsilon_n
    }}\Delta \tau_u 
  \end{equation*}
  Thus
  \begin{equation*}
  \sup_{r_i^n\le S_n}|\tau_{r_{i+1}^{n-}}-\tau_{r_i^n} | \le
  b \sup_{r_i^n\le S_n} |r_{i+1}^n-r_i^n| + 
  \sup_{r_i\le S_n}\sum_{\substack{
      r_i^n < u \le r_{i+1}^n \\
      \Delta \tau_u \le \epsilon_n
    }}\Delta \tau_u 
  \end{equation*}
  It has already been shown that the first term on the right converges
  to $0$ almost surely. Thus to prove \eqref{interJumpRanLim}, it
  suffices to prove that
  \begin{equation}\label{jumpSup}
    \lim_{n\to\infty} \sup_{r_i^n\le S_n}\sum_{\substack{
        r_i^n < u \le r_{i+1}^n \\
        \Delta \tau_u \le \epsilon_n
      }}\Delta \tau_u =0,
  \end{equation}
  almost surely, when $\epsilon_n\downarrow 0$ sufficiently
  quickly. Again, by Borel's lemma, \eqref{jumpSup} will follow if
  $\epsilon_n$ is chosen such that 
  \begin{equation}\label{jumpProbBound}
    \P\left(
      \sup_{r_i^n\le S_n}
      \sum_{\substack{
          r_i^n < u \le r_{i+1}^n \\
          \Delta \tau_u \le \epsilon_n
        }}\Delta \tau_u \ge \frac{1}{2^n}
    \right) < \frac{R}{2^n},
  \end{equation}
  for some $R>0$. 

  As before, suppress the superscripts on $r_i^n$ and the subscripts on
  $\epsilon_n$ and $S_n$. Recall that $r_{i+1}-r_i$ are exponential random
  variables with rate parameter $g(\epsilon)$. 
  Furthermore, the jump times of $\tau_s^n$ are independent of the
  small-jumps process
  \begin{equation*}
  \tau_s-\tau_s^n = \sum_{\substack{
      0\le r \le s \\
      \Delta\tau_r \le \epsilon
    }} \Delta\tau_r.
  \end{equation*}
  Define $h(\epsilon)$ by
  \begin{equation*}
  h(\epsilon) = \int_0^{\epsilon}t\lambda(dt).
  \end{equation*}
  Let $q(\epsilon)$ be the probability that
  $\sum_{\substack{
      r_i^n < u \le r_{i+1}^n \\
      \Delta \tau_u \le \epsilon_n
    }}\Delta \tau_u \ge \frac{1}{2^n}
  $. Define $\hat{q}(\epsilon)$ as the upper bound on $q(\epsilon)$ given
  by Markov's inequality:
  \begin{align}\label{qMarkov}
    \MoveEqLeft[0]
    q(\epsilon)=\P\left(
      \sum_{\substack{
          r_i^n < u \le r_{i+1}^n \\
          \Delta \tau_u \le \epsilon_n
        }}\Delta \tau_u \ge \frac{1}{2^n}
    \right)\le 
    \\ & \nonumber
    2^n \E\left[
      \sum_{\substack{
          r_i^n < u \le r_{i+1}^n \\
          \Delta \tau_u \le \epsilon_n
        }}\Delta \tau_u 
    \right] = 
    2^n\E\left[
      \int_{r_i}^{r_{i+1}}\int_0^{\epsilon}tN(dr,dt) 
    \right]
    =\frac{2^n
      h(\epsilon)
    }
    {
      g(\epsilon)
    }.
  \end{align}

  As in the proof of \eqref{interJumpDomBound}, the bound in \eqref{jumpProbBound}
  will be recast as a more tractable inequality. 

  Let $L$ be the geometric random variable defined by 
  \begin{equation*}
  L = \min\left\{i: \sum_{\substack{
        r_i < u \le r_{i+1} \\
        \Delta \tau_u \le \epsilon
      }}\Delta \tau_u \ge \frac{1}{2^n}
  \right\}.
  \end{equation*}
  So $L$ has probability given by $\P(L=k)= (1-q(\epsilon))^kq(\epsilon)
  $. As in the proof
  of \eqref{interJumpDomBound}, for any constant $M>0$, 
  \begin{align*}
    \MoveEqLeft
  \P\left(
    \sup_{r_i\le S}
    \sum_{\substack{
        r_i < u \le r_{i+1} \\
        \Delta \tau_u \le \epsilon
      }}\Delta \tau_u \ge \frac{1}{2^n}
  \right)
  = 
\\ &
\P(L\le K) \le \P(L\le M) + \P(M\le K). 
  \end{align*}
  The first term on the right can be bounded as 
  \begin{equation*}
  \P(L\le M)  = 1-(1-q(\epsilon))^{M+1}\le 1-(1-\hat{q}(\epsilon))^{M+1}.
  \end{equation*}
  Furthermore, as in the proof of \eqref{interJumpDomBound}, if 
  \begin{equation}\label{qBound}
    \left(2^nSg(\epsilon)+1\right)
    \log(1-\hat{q}(\epsilon)) < \log\left(1-\frac{1}{2^n}\right),
  \end{equation}
  the constant $M$ can be chosen such that 
  \begin{equation*}
  \P(L\le M) + \P(M\le K) \le \frac{2}{2^n}.
  \end{equation*}
  Thus, if \eqref{qBound} holds, then so does
  \eqref{jumpProbBound}. Note that $\hat{q}(\epsilon)\to 0$ as
  $\epsilon\to 0$, and thus $\log(1-\hat{q}(\epsilon))\to 0$ as
  well. Thus, for $\eqref{qBound}$ to hold for sufficiently small
  $\epsilon$, it suffices to show that
  $g(\epsilon)\log(1-\hat{q}(\epsilon))\to 0$ as $\epsilon\to 0$. 

  Using the power series expansion of $\log(1-\hat{q}(\epsilon))$
  implies that
  \begin{align*}
    \MoveEqLeft
    |g(\epsilon)\log(1-\hat{q}(\epsilon))| =
    g(\epsilon)\sum_{k=1}^{\infty}\frac{\hat{q}(\epsilon)^k}{k}
    =
    g(\epsilon)\sum_{k=1}^{\infty}\frac{2^{nk}h(\epsilon)^k}{g(\epsilon)^kk}
    \\&
    = 2^nh(\epsilon) \sum_{k=0}^{\infty}\frac{\hat{q}(\epsilon)^k}{k+1}
    \le 
    2^nh(\epsilon) \sum_{k=0}^{\infty}\hat{q}(\epsilon)^k
    =\frac{2^nh(\epsilon)}{1-\hat{q}(\epsilon)}.
  \end{align*}
  Now $\lim_{\epsilon\to 0}h(\epsilon) = 0$ implies that
  $\lim_{\epsilon\to 0} g(\epsilon)\log(1-\hat{q}(\epsilon)) = 0$. 
  Therefore \eqref{qBound} holds for sufficiently small $\epsilon$ and
  the proof is complete. 

  Now \eqref{finalRLim} and \eqref{finalTauLim} will be proved. As long
  as $S_n\to \infty$, the
  limit in \eqref{finalRLim} is immediate from \eqref{interJumpDomLim}
  since 
  \begin{equation*}
  S_n-\sup\{r_i^n:r_i^n \le S_n\} 
  \le \sup\{r_{i+1}^n-r_i^n:r_i^n\le S_n\}.
  \end{equation*}

  Now \eqref{finalTauLim} will be proved. If $b>0$, then
  \eqref{finalTauLim} follows for any sequence with $S_n\to\infty$. Thus, assume that
  $b=0$. Let $K_n = \max\{i:r_i^n \le S_n\}$ and assume that $\epsilon_{n-1}$ is fixed. The sequence in
  \eqref{finalTauLim} may be lower-bounded as
  \begin{gather*}
  \sup\{\tau_{r_i^n}^n:r_{i}^n\le S_n\}
  =\tau_{r_{K_n}^n} 
  \ge 
  \tau_{r_{K_n}^n}^n
  =\tau_{S_n}^n 
  \\
  \ge \tau_{S_n}^{n-1}\ge \epsilon_{n-1}N(S_n,(\epsilon_{n-1},\infty)).
  \end{gather*}
  The random measure term on the right is a Poisson process with rate
  $g(\epsilon_{n-1})$. Thus
  \begin{equation*}
  \P\left(
    \sup\{\tau_{r_i^n}^n:r_{i}^n\le S_n\} < n
  \right)
  \le 
  \P\left(
    \epsilon_{n-1}N(S_n,(\epsilon_{n-1},\infty)) < n
  \right)
  \end{equation*}
  Therefore, Borel's lemma implies that  it is sufficient to
  prove 
  \begin{equation}\label{PoissonBound}
    \P\left(
      \epsilon_{n-1}N(S_n,(\epsilon_{n-1},\infty)) < n
    \right) < \frac{1}{2^n}.
  \end{equation}

  The probability may be computed explicitly as
  \begin{align*}
    \MoveEqLeft
  \P\left(
    \epsilon_{n-1}N(S_n,(\epsilon_{n-1},\infty)) < n
  \right) = 
  \\ &
  e^{-S_ng(\epsilon_{n-1})}\sum_{0\le k <
    \frac{n}{\epsilon_{n-1}}}
  \frac{(S_ng(\epsilon_{n-1}))^k}{k}
  .
  \end{align*} 
  So, for fixed $\epsilon_{n-1}$, $S_n$ may be chosen sufficiently
  large so that \eqref{PoissonBound} holds, and the proof is complete. 
\end{IEEEproof}

\section{Proof of Lemma~\ref{lem:beta}} \label{pf:beta}
First note that  $\beta$ is analytic at $z$ if
$\int_0^{\infty}(e^{zt}-1)\lambda(dt)$ is. Furthermore, by the L\'evy-It\^o decomposition, 
\begin{equation*}
\E[e^{z\tau_s}] = e^{zbs}\E\left[
  \exp\left(z\int_0^{\infty}tN(s,dt)\right)
\right].
\end{equation*}
Thus, it suffices to prove the theorem for the case that $b=0$. 

Consider $z\in \dom(\beta)$, and let $r=\Re z$. It will be shown that
$\beta$ is analytic at $z$. Take any $\hat{r}\in (r,r_{\max})$ and any $y\in \C$ such that
$|y-z| < \hat{r}-r$. Provided that the sum of integrals below
converges absolutely, $\beta(y)$ can be derived from the power series
expansion around $z$:
\begin{align}
  \MoveEqLeft
  \nonumber
  \int_0^{\infty}\left(
    e^{zt}-1
  \right)\lambda(dt) + \sum_{k=1}^{\infty}\frac{1}{k!}(y-z)^k
  \int_0^{\infty}t^k e^{zt}\lambda(dt) 
  \\ \nonumber
  &= \int_0^{\infty}\left(
    e^{zt}-1
  \right)\lambda(dt) + \int_0^{\infty} \left(e^{(y-z)t}-1\right)
  e^{zt}\lambda(dt)
  \\
  \label{betaSeries}
  &=\beta(y).
\end{align}
Note that the form of the derivatives is immediate from the series
expansion.  

For absolute convergence, note that 
\begin{align}\label{absPowerSeries}
  \MoveEqLeft
  \int_0^{\infty}\left|
    e^{zt}-1
  \right|\lambda(dt) + \sum_{k=1}^{\infty}\frac{1}{k!}|y-z|^k
  \int_0^{\infty}\left|t^k e^{zt}\right|\lambda(dt) 
  =
  \\ & \nonumber
  \int_0^{\infty}\left|
    e^{zt}-1
  \right|\lambda(dt) + \sum_{k=1}^{\infty}\frac{1}{k!}|y-z|^k
  \int_0^{\infty}t^k e^{rt}\lambda(dt). 
\end{align}

The first term on the right is seen to be finite as follows. If $t\in
[0,1]$, then
\begin{equation*}
\left|
  e^{zt}-1
\right| = \left|
  \sum_{k=1}^{\infty}\frac{(zt)^k}{k!}
\right|\le t \sum_{k=1}^{\infty}\frac{|z|^k}{k!}
=t\left(e^{|z|}-1\right). 
\end{equation*} 
On the other hand, the triangle inequality implies that 
\begin{equation*}
\left|
  e^{zt}-1
\right| \le \left|e^{zt}\right| + 1 = e^{rt} + 1.
\end{equation*}
Therefore, the first integral on the right of \eqref{absPowerSeries} is bounded as 
\begin{align}\label{diffIntegralBound}
  \MoveEqLeft
  \int_0^{\infty} \left|
    e^{zt}-1
  \right|\lambda(dt) 
  \\ & \nonumber
\le \left(e^{|z|}-1\right)\int_0^1 t \lambda(dt)
  + \int_1^{\infty} \left(e^{rt} + 1\right)\lambda(dt)<\infty,
\end{align}
where the last inequality follows since $r<r_{\max}$. 

Now the integral in the sum on the right of \eqref{absPowerSeries}
will be bounded. First note that the integrand is bounded as 
\begin{equation}\label{integrandBound}
  t^ke^{rt}=t^k
  e^{-(\hat{r}-r)t}e^{\hat{r} t} \le \left(
    \frac{k}{e(\hat{r}-r)}
  \right)^k e^{\hat{r}t},
\end{equation}
where the inequality follows from maximizing
$te^{-(\hat{r}-r)t}$. Thus, the integral is bounded as 
\begin{equation}\label{integralBound}
  \int_{0}^{\infty}t^ke^{rt}\lambda(dt) 
  \le e^r \int_0^1 t\lambda(dt) + \left(
    \frac{k}{e(\hat{r}-r)}
  \right)^k\int_1^{\infty}e^{\hat{r}t}\lambda(dt)
\end{equation}
Thus, to prove that the sum on the right of \eqref{absPowerSeries}
converges, it suffices to prove that 
\begin{equation*}
\sum_{k=1}^{\infty} \frac{|y-z|^k}{k!}\left(
  \frac{k}{e(\hat{r}-r)}  
\right)^k<\infty
\end{equation*}
Now, the Stirling approximation bound, $k! \ge \sqrt{2\pi k} \left(
  \frac{k}{e}
\right)^k$, shows that 
\begin{equation*}
\frac{1}{k!} \left(
  \frac{k}{e(\hat{r}-r)}  
\right)^k
\le \frac{1}{\sqrt{2\pi k} (\hat{r}-r)^k}
\le \frac{1}{(\hat{r}-r)^k}. 
\end{equation*}
Since $|y-z| < \hat{r}-r$, the bound follows as
\begin{equation*}
\sum_{k=1}^{\infty} \frac{|y-z|^k}{k!}\left(
  \frac{k}{e(\hat{r}-r)}  
\right)^k
\le \sum_{k=1}^{\infty} \left(\frac{|y-z|}{\hat{r}-r}\right)^k <
\infty. 
\end{equation*}
Thus, the power series expansion, \eqref{betaSeries} holds, and
$\beta$ is analytic at $z$.


Now \eqref{betaExp} will be proved. 
The proof is similar to the proof of Theorem 2.3.8 in \cite{applebaumlevy2004}. 

First, the function $t$ will be
approximated by step functions over $(0,\infty)$. The construction is
similar to the approach in the proof of Theorem 1.17 in \cite{rudinreal1987}. Consider a sequence
$\gamma_n\downarrow 0$ at a rate to be specified later. Let $k_n(t)$
be the unique integer such that $k \gamma_n \le t < (k+1)
\gamma_n$. Define the function $\varphi_n(t)$ by 
\begin{equation*}
\varphi_n(t) = \left\{
  \begin{matrix}
    k_n(t) \gamma_n & \qquad t\in (0,n) \\
    n & \qquad t \ge n.
  \end{matrix}
\right.
\end{equation*}
Then $\varphi_n(t)$ is a simple function such that $\varphi_n(t) =0$
for $t\in (0,\gamma_n)$, $t-\gamma_n< \varphi_n(t) \le t$ for $t\in
[\gamma_n,n]$, and $\varphi_n(t)\le t$ for $t>0$. The formula, \eqref{betaExp}, is a consequence of the following chain of equalities
\begin{align}
  \MoveEqLeft \nonumber
  \E\left[
    \exp\left(z\int_0^{\infty}tN(s,dt)\right)
  \right] 
\\ \label{betaExpApprox1}
  &= \lim_{n\to\infty} \E\left[
    \exp\left(z\int_0^{\infty} \varphi_n(t) N(s,dt)
    \right)
  \right] \\
  \label{betaExpApprox2}
  &= \lim_{n\to\infty} \exp\left(
    s\int_0^{\infty}\left(
      e^{z\varphi_n(t)} -1 
    \right)\lambda(dt)
  \right)
  \\&=
  \label{betaExpApprox3}
  \exp\left(
    s\int_0^{\infty}\left(
      e^{z t} -1 
    \right)\lambda(dt)
  \right).
\end{align}

The first equation is the most challenging, and will be handled last. To prove \eqref{betaExpApprox2}, note that $z\varphi_n(t)$ is a simple function. Thus, there are constants $c_i\in \C$ and disjoint $\lambda$-measurable sets, $A_i$, such that 
\begin{equation*}
z\varphi_n(t) = \sum_{i=1}^q c_i \ind_{A_i}(t). 
\end{equation*}

Since $\varphi_n(t)=0$ over $(0,\gamma_n)$, it follows that $0$ is not in the closure of any $A_i$. Thus, the integral on the right of \eqref{betaExpApprox1} may be written as 
\begin{equation*}
\int_0^{\infty}z\varphi_n(t)N(s,dt) = \sum_{i=1}^q c_i N(s,A_i),
\end{equation*}
where $N(s,A_i)$ are independent Poisson processes with rate $\lambda(A_i)$. Thus, the expectation on the right of \eqref{betaExpApprox1} may be calculated as 
\begin{align*}
  \MoveEqLeft
  \E\left[
    \exp\left(z\int_0^{\infty} \varphi_n(t) N(s,dt)
    \right)
  \right] 
  \\
&= \prod_{i=1}^q \E\left[
    \exp(c_i N(s,A_i))
  \right] \\
  &= \prod_{i=1}^q  \exp(-s\lambda(A_i))\sum_{k=0}^{\infty}
  \frac{(s\lambda(A_i))^k}{k!} e^{c_i k}\\
  &= \prod_{i=1}^q \exp\left(s\lambda(A_i) \left(
      e^{c_i} - 1
    \right)\right)
  \\
  &= 
  \exp\left(
    s\sum_{i=1}^q \left(e^{c_i}-1\right) \lambda(A_i)
  \right)
  \\
  &= \exp\left(s\int_0^{\infty} \left(
      e^{z\varphi_n(t)} - 1
    \right)\lambda(dt)
  \right).
\end{align*}
Thus, \eqref{betaExpApprox2} holds.  

To prove \eqref{betaExpApprox3}, note that \eqref{diffIntegralBound}
implies that $\left|e^{z\varphi_n(t)}-1\right|$ is bounded above by a
function with a finite integral. Thus, Lebesgue's dominated
convergence theorem implies that 
\begin{equation*}
\lim_{n\to\infty} \int_0^{\infty} \left(
  e^{z\varphi_n(t)} -1
\right)\lambda(dt) = \int_0^{\infty}\left( e^{zt}-1\right)\lambda(dt)
\end{equation*}  
and so \eqref{betaExpApprox3} holds. 

Finally, \eqref{betaExpApprox1} must be proved. First, it will be
shown that 
\begin{equation}\label{simpleProcLim}
  \lim_{n\to\infty} \int_0^{\infty} \varphi_n(t) N(s,dt) = \int_0^t t
  N(s,dt),
  \qquad a.s. 
\end{equation}
Then, dominated convergence will be
applied. 

Assume that $\gamma_{n-1}$ is fixed. The difference of the right and
left of \eqref{simpleProcLim} may be bounded as 
\begin{align}\label{ProcDiffBound}
  \MoveEqLeft
  0\le \int_0^{\infty}(t-\varphi_n(t)) N(s,dt)\le \\ & \nonumber 
  \int_0^{\gamma_{n-1}} t N(s,dt) + \gamma_n \int_{\gamma_{n-1}}^n
  N(s,dt) + \int_{n}^{\infty}(t-n)N(s,dt). 
\end{align}
To bound the first term on the right of \eqref{ProcDiffBound}, note
that 
\begin{equation*}
\int_0^{\gamma_{1}} t N(s,dt) = \sum_{i=1}^{\infty}
\int_{\gamma_{i+1}}^{\gamma_i} t N(s,dt) <\infty, \qquad
\textrm{almost surely}. 
\end{equation*}
Thus, the the first term on the right of \eqref{ProcDiffBound} may be
expressed as the tail sum:
\begin{equation*}
\int_0^{\gamma_{n-1}} t N(s,dt) = \sum_{i=n-1}^{\infty}
\int_{\gamma_{i+1}}^{\gamma_i} t N(s,dt)
\end{equation*}
which converges to $0$ almost surely, provided that $\gamma_n\downarrow 0$ sufficiently quickly. (See \cite{applebaumlevy2004}.)  

Now consider the second term on the right of \eqref{ProcDiffBound}. For fixed $\gamma_{n-1}$, the next term
$\gamma_n$ may be chosen sufficiently small to give the following
probability bound:
\begin{align*}
 \MoveEqLeft
\P\left(
  \gamma_{n}\int_{\gamma_{n-1}}^n N(s,dt) \ge 2^{-n}
\right) = 
\\ &
e^{-s \lambda([\gamma_{n-1},n))}\sum_{k\ge
  \frac{1}{\gamma_n 2^n}} 
\frac{
  (s \lambda([\gamma_{n-1},n)))^k
}{k!} < \frac{1}{2^n}.
\end{align*}
Thus, by Borel's lemma, the second term converges to $0$ almost
surely. 

The last term on the right of \eqref{ProcDiffBound} is $0$ if $\tau_s
<n$, which holds for sufficiently large $n$ almost surely. Thus
\eqref{simpleProcLim} holds. 

Now it will be shown that Lebesgue's dominated convergence applies to
\eqref{betaExpApprox1}. Note that the function on the right has magnitude given by
\begin{align*}
\left|\exp\left(z\int_0^{\infty} \varphi_n(t) N(s,dt)
  \right)\right| &= \exp\left(r\int_0^{\infty} \varphi_n(t) N(s,dt)
\right) .
\end{align*}
Thus, it suffices to show that the term on the right has finite
expectation. If $r\le 0$, then the term is bounded above by $1$ and so
finiteness is immediate. So, consider the case that $r>0$. Here the magnitude is bounded above by 
\[
\exp\left(r \int_0^{\infty} t N(s,dt)\right) = \exp(r\tau_s).
\]
The expectation of this term may be evaluated, as long as the following
equalities can be proved:
\begin{align}
  \nonumber
  \E\left[
    \exp\left(r \tau_s\right)
  \right] &= \E\left[
    \sum_{k=0}^{\infty}\frac{(r\tau_s)^k}{k!} 
  \right] \\
  \label{expExp1}
  &= \sum_{k=0}^{\infty} \frac{r^k}{k!} \E\left[\tau_s^k\right] \\
  \label{expExp2}
  &=1+ s\sum_{k=1}^{\infty} \frac{r^k}{k!} \int_0^{\infty}t^k \lambda(dt)
  \\
  \label{expExp3}
  &=1+ s \int_0^{\infty}\left(e^{rt}-1\right)\lambda(dt)
  \\
  \nonumber
  & < \infty
\end{align}
The inequality follows since $r<r_{\max}$, while the first
equality is just the definition of the exponential function. 
Note that $0< r \in\dom(\beta)$
implies that $0\in\dom(\beta)$. Thus, $\beta$ is analytic at $0$, and
the argument above implies that the integrals on the right of
\eqref{expExp2} are finite. 
Therefore, \eqref{expExp3} follows from non-negativity and Fubini's
theorem. Furthermore, provided that \eqref{expExp2} holds, \eqref{expExp1}
will hold by Fubini's theorem. 

Now, the only remaining equality, \eqref{expExp2}, will be shown.  Recall the definition of $\tau_s^n$
from \eqref{tauNdef}, and recall that in this case, $b=0$. By
construction $\tau_s^n\uparrow \tau_s$ almost surely. Using Lebesgue's dominated
convergence theorem and then Theorem
2.3.8 of \cite{applebaumlevy2004}, the
following equalities hold for $k\ge 1$:
\begin{equation*}
\E\left[\tau_s^k\right] = \lim_{n\to\infty} \E\left[(\tau_s^n)^k\right] =
s\lim_{n\to\infty} \int_{\epsilon_n}^{\infty} t^k \lambda(dt) 
= s\int_0^{\infty} t^k \lambda(dt). 
\end{equation*}
Thus, \eqref{expExp2} has been shown, and so \eqref{betaExp} holds.

For the matrix case, the following generalization of \eqref{betaExp} is helpful:
\begin{equation}\label{betaExpDeriv}
  \E\left[\tau_s^ke^{z\tau_s}\right] = \frac{\partial^k }{\partial z^k} e^{s\beta(z)}.
\end{equation}
It is proved using Cauchy's integral formula:
\begin{align*}
  \E\left[
    \tau_s^k e^{z\tau_s} 
  \right]&= 
  \E
  \left[
    \frac{k!}{2\pi i} \oint_C \frac{e^{y\tau_s}}{(y-z)^{k+1}}dy
  \right]
  \\
  &= \frac{k!}{2\pi i} \oint_C
  \frac{\E\left[e^{y\tau_s}\right]}{(y-z)^{k+1}}dy \\
  &= \frac{k!}{2\pi i} \oint_C
  \frac{e^{s\beta(y)}}{(y-z)^{k+1}}dy
  \\
  &= \frac{\partial^k}{\partial z^k} e^{s\beta(z)}.
\end{align*}
The only equality requiring justification is the second, which follows
from Fubini's theorem provided that $\Re y < \hat{r}$ for all $y$ on
the contour.

Now the definition of $\beta$ and \eqref{betaExp} will be extended to
matrices. For full generality, the $b$ term will be included. Consider
a Jordan decomposition, $A = V^{-1}JV$, and let $J_i$ be an $m\times m$ Jordan
block corresponding to eigenvalue $\rho$. For $J_i$,
\eqref{matrixBeta} may be evaluated as 
\begin{equation*}
\beta(J_i) = bJ_i + \int_0^{\infty}
\left(
  e^{\rho t}
  \begin{bmatrix}
    1 & t & \cdots & \frac{t^{m-1}}{(m-1)!} \\
    & \ddots & \ddots & \vdots \\
    && 1 & t \\
    &&& 1
  \end{bmatrix}
  -I
\right)\lambda(dt).
\end{equation*}
Since $\rho\in\spec(A)\subset \dom(\beta)$, it follows that the
integral converges for all entries. Therefore, $\beta(A)$ may be
computed as 
\begin{align*}
\MoveEqLeft
\beta(A) = bA + \int_0^{\infty}\left(e^{At}-I\right) \lambda(dt) =
\\ &
V^{-1}
\left(bJ+\int_0^{\infty}\left(e^{Jt}-I\right)\lambda(dt)\right)V
=V^{-1} \beta(J)V.
\end{align*}


The proof of
\eqref{betaExp} generalizes in a straightforward manner when $z$ is replaced by the Jordan block, $J_i$. Again, assume that $b=0$.  The following equalities must be shown
\begin{align}\label{matrixChain1}
  \E\left[e^{J_i \tau_s}\right] &= \lim_{n\to\infty}\E\left[\exp\left(
      J_i\int_0^{\infty}\varphi_n(t)N(s,dt)
    \right)
  \right] 
  \\ & \label{matrixChain2}
  = \lim_{n\to\infty} \exp\left(
    s\int_0^{\infty} \left(e^{J_i\varphi_n(t)}-I\right)\lambda(dt)
  \right)
  \\ & \label{matrixChain3}
  =e^{s\beta(J_i)}.
\end{align}
If \eqref{matrixChain1}-\eqref{matrixChain3} hold, then
then $\E[e^{A\tau_s}] = V^{-1}\E[e^{J\tau_s}] V = V^{-1}e^{s \beta(J)}V = e^{s\beta(A)}$.

The second equality, \eqref{matrixChain2}, holds by formally following the steps in the derivation of \eqref{betaExpApprox2}. The first and third equalities will hold as long as the off-diagonal terms may be bounded in order to apply Lebesgue's dominated convergence theorem.


Consider  \eqref{matrixChain1}. Let $r=\Re \rho$ and let $k$ be a positive integer. If $r\ge 0$ note that
\begin{align*}
  \MoveEqLeft
  \left|
    \left(\int_0^{\infty} \varphi_n(t)N(s,dt)\right)^k \exp\left(\rho \int_0^{\infty}\varphi_n(t)N(s,dt)\right)
  \right|\le 
  \\ &
  \left(\int_0^{\infty}tN(s,dt)\right)^k \exp\left(r \int_0^{\infty}tN(s,dt)\right),
\end{align*}
while if $r<0$, the left is bounded by a constant, as in \eqref{integrandBound}. In either case, the upper bound has finite expectation, according to \eqref{betaExpDeriv}. Thus, the \eqref{matrixChain1} must hold. 

Now the third equality, \eqref{matrixChain3}, will be proved. Note that 
\begin{equation*}
\left|\varphi_n(t)^k e^{\rho \varphi_n(t)}\right| \le t^k e^{rt}K,
\end{equation*}
where $K\ge \max\{e^{-r\gamma_n},1\}$. The upper bound has a finite $\lambda$-integral  for $k\ge 1$ because of \eqref{integralBound}. 
Thus \eqref{matrixChain3} holds, and the proof is complete. \hfill\IEEEQED

\section{Proof of Lemma~\ref{lem:mappings}}\label{pf:mappings}

Define the matrix $Z$ by 
\begin{equation*}
Z = \begin{bmatrix} I \\ 0 \end{bmatrix} P \begin{bmatrix}I
  &0 \end{bmatrix} 
\end{equation*}
and define the matrix $\tilde{A}$ by
\begin{equation*}
\tilde{A} = \begin{bmatrix} A & I \\ 0 & 0 \end{bmatrix}.
\end{equation*}
Note that $e^{\tilde{A} t}$ is given by 
\begin{equation*}
e^{\tilde{A} t} = \begin{bmatrix}
  e^{At} & \int_0^t e^{Ar}dr \\
  0 & I
\end{bmatrix}.
\end{equation*}

Thus, the matrix-valued mappings may be written as 
\begin{align*}
\MoveEqLeft
\begin{bmatrix}
  F(P) & G(P) \\
  G(P)^\tp & H(P)
\end{bmatrix}
=
\\ &
b\left(
  \tilde{A}^\tp Z + Z \tilde{A}
\right)+
\int_0^{\infty}
\left(e^{\tilde{A}^\tp t} Z e^{\tilde{A}t}
  -Z
\right)
\lambda(dt)
\end{align*}
Since $e^{\tilde{A}^\tp t}\otimes e^{\tilde{A}^\tp t} =
e^{\tilde{A}^\tp \oplus \tilde{A}^\tp t}$, the equation may be
vectorized as 
\begin{align*}
\MoveEqLeft
  \vec\left(
    \begin{bmatrix}
      F(P) & G(P) \\
      G(P)^\tp & H(P)
    \end{bmatrix}
  \right)  
\\ 
&= 
  \left(b \tilde{A}^\tp \oplus\tilde{A}^\tp + \int_0^{\infty}
    \left(
      e^{\tilde{A}^\tp \oplus \tilde{A}^\tp t} - I
    \right)\lambda(dt)\right) \vec(Z) 
  \\ 
  &= \beta(\tilde{A}^\tp
  \oplus\tilde{A}^\tp ) \vec(Z).
\end{align*}
Thus according to Lemma~\ref{lem:beta}, $F$, $G$, and $H$ are
well defined, as long as $\spec(\tilde{A}^{\tp}\oplus
\tilde{A}^{\tp})\subset \dom(\beta)$. By construction, the spectrum is
given by
\begin{align*}
\spec(\tilde{A}^{\tp}\oplus
\tilde{A}^{\tp}) &= \spec(\tilde{A}^\tp) + \spec(\tilde{A}^\tp) \\
&= \{0\}\cup\spec(A)\cup(\spec(A)+\spec(A)).
\end{align*}

Let $r=\max\{\Re \mu: \mu\in\spec(A)\}$. If $r\le 0$, then the maximum
real part of any eigenvalue of $\tilde{A}^{\tp}\oplus
\tilde{A}^{\tp}$ is $0$. If $r>0$, then the corresponding maximum real
part must be $2r$. Since $\{0\}\cup\spec(2A)\subset \dom(\beta)$, it
follows that $\spec(\tilde{A}^{\tp}\oplus
\tilde{A}^{\tp})\subset \dom(\beta)$, and so the mappings are defined.

Furthermore,  the relevant expectations may be vectorized and
evaluated using \eqref{betaM}:
\begin{align*}
\MoveEqLeft
\vec\left(
  \E\left[
    e^{\tilde{A}^\tp \tau_s} Z e^{\tilde{A}^\tp \tau_s}
  \right]
\right)
\\
&= \E\left[
  e^{\tilde{A}^\tp \oplus \tilde{A}^\tp \tau_s }
\right]\vec(Z) 
\\
&= \vec(Z) + s\beta(\tilde{A}^\tp \oplus \tilde{A}^\tp)\vec(Z)
+ O(s^2).
\end{align*}

The proof for $g$ is similar, noting that 
\begin{equation*}
\vec\left(\int_0^t e^{At}MM^\tp e^{A^\tp t} dt\right)
= \begin{bmatrix}I & 0 \end{bmatrix}
\left(
  e^{\hat{A}t} - I
\right) \begin{bmatrix} 0 \\ I \end{bmatrix} \vec(MM^\tp),
\end{equation*}
where 
\begin{equation*}
\hat{A} = \begin{bmatrix}
  A\oplus A & I \\
  0 & 0
\end{bmatrix}.
\end{equation*}
\hfill\IEEEQED
\bibliography{bibLoc}

\begin{thebibliography}{10}
\providecommand{\url}[1]{#1}
\csname url@rmstyle\endcsname
\providecommand{\newblock}{\relax}
\providecommand{\bibinfo}[2]{#2}
\providecommand\BIBentrySTDinterwordspacing{\spaceskip=0pt\relax}
\providecommand\BIBentryALTinterwordstretchfactor{4}
\providecommand\BIBentryALTinterwordspacing{\spaceskip=\fontdimen2\font plus
\BIBentryALTinterwordstretchfactor\fontdimen3\font minus
  \fontdimen4\font\relax}
\providecommand\BIBforeignlanguage[2]{{%
\expandafter\ifx\csname l@#1\endcsname\relax
\typeout{** WARNING: IEEEtran.bst: No hyphenation pattern has been}%
\typeout{** loaded for the language `#1'. Using the pattern for}%
\typeout{** the default language instead.}%
\else
\language=\csname l@#1\endcsname
\fi
#2}}

\bibitem{lavalletime2007}
S.~M. LaValle and M.~B. Egerstedt, ``On time: Clocks, chronometers, and
  open-loop control,'' in \emph{IEEE Conference on Decision and Control}, 2007.

\bibitem{carverstateestimation2013}
S.~G. Carver, E.~S. Fortune, and N.~J. Cowan, ``State-estimation and
  cooperative control with uncertain time,'' in \emph{American Control
  Conference}, 2013.

\bibitem{lamperskitimechanged2013}
A.~Lamperski and N.~J. Cowan, ``Time-changed linear quadratic regulators,'' in
  \emph{European Control Conference}, 2013.

\bibitem{veraarttime2010}
A.~E.~D. Veraart and M.~Winkel, ``Time change,'' in \emph{Encyclopedia of
  Quantitative Finance}.\hskip 1em plus 0.5em minus 0.4em\relax Wiley, 2010,
  vol.~4, pp. 1812--1816.

\bibitem{clarksubordinated1973}
P.~K. Clark, ``A subordinated stochastic process model with finite variance for
  speculative prices,'' \emph{Econometrica}, vol.~41, pp. 135--155, 1973.

\bibitem{aneorder2000}
T.~An\'e and H.~Geman, ``Order flow, transaction clock, and normality of asset
  returns,'' \emph{The Journal of Finance}, vol.~55, no.~5, pp. 2259--2284,
  2000.

\bibitem{carrtimechanged2004}
P.~Carr and L.~Wu, ``Time-changed {L}\'evy processes and option pricing,''
  \emph{Journal of Financial Economics}, vol.~71, pp. 113--141, 2004.

\bibitem{eaglemanhuman2008}
D.~M. Eagleman, ``Human time perception and its illusions,'' \emph{Current
  Opinion in Neurobiology}, vol.~18, no.~2, 2008.

\bibitem{jazayeritemporal2010}
M.~Jazayeri and M.~N. Shadlen, ``Temporal context calibrates interval timing,''
  \emph{Nature Neuroscience}, vol.~13, no.~8, 2010.

\bibitem{hudsonoptimal2008}
T.~E. Hudson, L.~T. Maloney, and M.~S. Landy, ``Optimal compensation for
  temporal uncertainty in movement planning,'' \emph{{PLoS} Computational
  Biology}, vol.~4, no.~7, 2008.

\bibitem{kushnerstability1969}
H.~J. Kushner and L.~Tobias, ``On the stability of randomly sampled systems,''
  \emph{IEEE Transactions on Automatic Control}, vol.~14, no.~4, pp. 319--324,
  1969.

\bibitem{wittenmarktiming1995}
B.~Wittenmark, J.~Nilsson, and M.~T\"orngren, ``Timing problems in real-time
  control systems,'' in \emph{American Control Conference}, 1995.

\bibitem{skafanalysis2009}
J.~Skaf and S.~Boyd, ``Analysis and synthesis of state-feedback controllers
  with timing jitter,'' \emph{IEEE Transactions on Automatic Control}, vol.~54,
  no.~3, pp. 652--657, 2009.

\bibitem{hespanhasurvey2007}
J.~P. Hespanha, P.~Naghshtabrizi, and Y.~Xu, ``A survey of recent results in
  networked control systems,'' \emph{Proceedings of the IEEE}, vol.~95, no.~1,
  pp. 138--162, 2007.

\bibitem{adesstochastic2000}
M.~Ad\`es, P.~E. Caines, and R.~P. Malham\`e, ``Stochastic optimal control
  under poisson-distributed observations,'' \emph{IEEE Transactions on
  Automatic Control}, vol.~45, no.~1, pp. 3--13, 2000.

\bibitem{frerisfundamental2007}
N.~M. Freris and P.~R. Kumar, ``Fundamental limits on synchronization of affine
  clocks in networks,'' in \emph{IEEE Conference on Decision and Control},
  2007.

\bibitem{simeonedistributed2008}
O.~Simeone, U.~Spagnolini, Y.~Bar-Ness, and S.~H. Strogatz, ``Distributed
  synchronization in wireless networks,'' \emph{IEEE Signal Processing
  Magazine}, vol.~25, no.~5, pp. 81--97, 2008.

\bibitem{lorandstability2003}
C.~Lorand and P.~H. Bauer, ``Stability analysis of closed-loop discrete-time
  systems with clock frequency drifts,'' in \emph{American Control Conference},
  2003.

\bibitem{singhlqr2011}
R.~Singh and V.~Gupta, ``On {LQR} control with asynchronous clocks,'' in
  \emph{IEEE Conference on Decision and Control and European Control
  Conference}, 2011.

\bibitem{applebaumlevy2004}
D.~Applebaum, \emph{L\'evy Processes and Stochastic Calculus}.\hskip 1em plus
  0.5em minus 0.4em\relax Cambridge University Press, 2004.

\bibitem{highamfunctions2008}
N.~J. Higham, \emph{Functions of Matrices: Theory and Computation}.\hskip 1em
  plus 0.5em minus 0.4em\relax SIAM, 2008.

\bibitem{cvitanicoptimal2008}
J.~Cvitani\'c, V.~Polimenis, and F.~Zapatero, ``Optimal portfolio allocation
  with higher moments,'' \emph{Annals of Finance}, vol.~4, pp. 1--28, 2008.

\bibitem{madanvariance1998}
D.~B. Madan, P.~P. Carr, and E.~C. Chang, ``The variance gamma process and
  option pricing,'' \emph{European Finance Review}, vol.~2, pp. 79--105, 1998.

\bibitem{barndorffnielsenprocesses1998}
O.~E. Barndorff-Nielsen, ``Processes of normal inverse {G}aussian type,''
  \emph{Finance and Stochastics}, vol.~2, pp. 41--68, 1998.

\bibitem{michaelgenerating1976}
J.~R. Michael, W.~R. Schucany, and R.~W. Haas, ``Generating random variates
  using transformations with multiple roots,'' \emph{The American
  Statistician}, vol.~30, no.~2, pp. 88--90, 1976.

\bibitem{satolevy1999}
K.~Sato, \emph{L\'evy Processes and Infinitely Divisible Distributions}.\hskip
  1em plus 0.5em minus 0.4em\relax Cambridge University Press, 1999.

\bibitem{flemingcontrolled2006}
W.~H. Fleming and H.~M. Soner, \emph{Controlled Markov Processes and Viscosity
  Solutions}, 2nd~ed.\hskip 1em plus 0.5em minus 0.4em\relax Springer, 2006.

\bibitem{protterstochastic2004}
P.~E. Protter, \emph{Stochastic Integration and Differential Equations},
  2nd~ed.\hskip 1em plus 0.5em minus 0.4em\relax Springer, 2004.

\bibitem{todorovoptimal2002}
E.~Todorov and M.~I. Jordan, ``Optimal feedback control as a theory of motor
  coordination,'' \emph{Nature Neuroscience}, vol.~5, no.~11, 2002.

\bibitem{harrissignal-dependent1998}
C.~Harris and D.~Wolpert, ``Signal-dependent noise determines motor planning,''
  \emph{Nature}, vol. 394, no. 6695, pp. 780--784, 1998.

\bibitem{rogersoptimal2013}
L.~C.~G. Rogers, \emph{Optimal Investment}.\hskip 1em plus 0.5em minus
  0.4em\relax Springer, 2013.

\bibitem{jadbabaiecoordination2003}
A.~Jadbabaie, J.~Lin, and A.~S. Morse, ``Coordination of groups of mobile
  autonomous agents using nearest neighbor rules,'' \emph{IEEE Transactions on
  Automatic Control}, vol.~48, pp. 988--1001, 2003.

\bibitem{olfatisaberconsensus2004}
R.~Olfati-Saber and R.~M. Murray, ``Consensus problems in networks of agents
  with switching topology and time-delays,'' \emph{IEEE Transactions on
  Automatic Control}, vol.~49, pp. 1520--1533, 2004.

\bibitem{gemanstochastic2002}
H.~Geman, D.~B. Madan, and M.~Yor, ``Stochastic volatility, jumps and hidden
  time changes,'' \emph{Finance and Stochastics}, vol.~6, pp. 63--90, 2002.

\bibitem{shibayesian2013}
Z.~Shi, R.~M. Church, and W.~H. Meck, ``Bayesian optimization of time
  perception,'' \emph{Trends in Cognitive Sciences}, vol.~17, no.~11, pp.
  556--564, 2013.

\bibitem{ankaralihaptic2014}
M.~M. Ankaral{\i}, H.~T. \d{S}en, A.~De, A.~M. Okamura, and N.~J. Cowan,
  ``Haptic feedback enhances rhythmic motor control performance by reducing
  variability, not convergence time,'' \emph{Journal of Neurophysiology}, 2014,
  in press.

\bibitem{rudinreal1987}
W.~Rudin, \emph{Real and Complex Analysis}.\hskip 1em plus 0.5em minus
  0.4em\relax WCB/McGraw-Hill, 1987.

\end{thebibliography}
\end{document}